\documentclass[a4paper,12pt]{amsart}\usepackage{amsfonts,amssymb,bbm}
\usepackage{enumerate,mathrsfs,xfrac,xspace} \UseRawInputEncoding
\usepackage[bookmarksdepth=3]{hyperref}

\def\C{\Bbb{C}}
\def\k{\Bbbk}
\def\K{\mathbbm{K}}
\def\bk{{\bar{\k}}}

\def\N{\Bbb{N}}\def\P{\Bbb{P}}\def\Q{\Bbb{Q}}\def\R{\Bbb{R}}\def\Z{\Bbb{Z}}

\def\di{\partial}

\def\bl{\langle}\def\br{\rangle}

\def\liml{\lim\limits}

\newcommand{\quot}[2]{{\footnotesize\left.\raisebox{1.2ex}{$#1$}\!\! \ensuremath\diagup \!\!\raisebox{-1.2ex}{$#2$}\right.}}
\newcommand{\quots}[2]{{\footnotesize\left.\raisebox{0.4ex}{$#1$}\! / \!\raisebox{-0.4ex}{$#2$}\right.}}

\renewcommand{\stackrel}[2]{\ \!\lower 0.3ex \hbox{$\mathrel{\mathop{#2}\limits^{#1}}$}\ \!}

\def\tf{{\tilde{f}}}

\def\tg{\tilde{g}}\def\tM{\tilde{M}}

\def\txi{{\tilde{\xi}}}

\def\hotimes{\hat{\otimes}}

\def\hg{{\hat{g}}}
\def\hk{{\hat\k}}\def\hR{{\widehat{R}}}

\def\al{\alpha}\def\be{\beta}
\def\ze{\zeta}
\def\la{\lambda}

\def\cA{\mathscr A}\def\ca{\mathfrak a}
\def\cC{{\mathscr C}}

\def\cG{\mathscr G}
\def\cK{{\mathscr K}\!}\def\cL{{\mathscr L}}\def\cP{\frak P}\def\cR{\mathscr{R}}

\def\cm{{\frak m}}

\def\one{{1\hspace{-0.1cm}\rm I}}

\newcommand{\ber}{\begin{array}{l}}\newcommand{\eer}{\end{array}}
\newcommand{\bpm}{\begin{pmatrix}}\newcommand{\epm}{\end{pmatrix}}
\newcommand{\bbm}{\begin{bmatrix}}\newcommand{\ebm}{\end{bmatrix}}
\newcommand{\bM}{\begin{matrix}}\newcommand{\eM}{\end{matrix}}
\newcommand{\bee}{\begin{enumerate}}\newcommand{\eee}{\end{enumerate}}
\newcommand{\bei}{\begin{itemize}}\newcommand{\eei}{\end{itemize}}

\def\sset{\!\subset\!}\def\sseteq{\!\subseteq\!}\def\ssetneq{\!\subsetneq\!}
\def\spset{\!\supset\!}\def\spseteq{\!\supseteq\!}
\def\smin{\!\setminus\!}\def\={\!=\!}
\def\scdot{\!\cdot\!}
\def\sin{\!\in\!}

\def\wrt{with respect to }\def\iff{if and only if }

\def\Maps{{\rm Maps}\left(X,Y\right)}\def\MapsK{{\rm Maps}\left(X_\K,Y_\K\right)}\def\MapX{{\rm Maps}\left(X,\!(\k^m,o)\right)}
\def\Mapk{{\rm Maps}\left((\k^n,o),\!(\k^m,o)\right)}\def\MapR{{\rm Maps}\!\left((\R^n,o),\!(\R^m,o)\right)\!}\def\MapC{{\rm Maps}\!\left((\C^n,o),\!(\C^m,o)\right)\!}

\def\Spec{\rm{Spec}}

\def\RmX{R^{\oplus m}_X}\def\RrX{R^{\oplus r}_X}
\def\RrX{R^{\oplus r}_X}\def\RmXK{R^{\oplus m}_{X,\K}}

\newcommand{\beq}{\begin{equation}}
\newcommand{\eeq}{\end{equation}}
\numberwithin{equation}{section}

\newtheorem{Lemma}{Lemma}[section]\newcommand{\bel}{\begin{Lemma}}\newcommand{\eel}{\end{Lemma}}
\newtheorem{Example}[Lemma]{Example}\newcommand{\bex}{\begin{Example}\rm}
\newcommand{\eex}{\end{Example}}
\newtheorem{Proposition}[Lemma]{Proposition}\newcommand{\bprop}{\begin{Proposition}}\newcommand{\eprop}{\end{Proposition}}
\newtheorem{Property}[Lemma]{Property}\newcommand{\bproperty}{\begin{Property}}\newcommand{\eproperty}{\end{Property}}
\newtheorem{Definition-Proposition}[Lemma]{Definition-Proposition}

\def\bpr{~\\{\em Proof.\ }}
\newcommand{\epr}{{\hfill\ensuremath\blacksquare}\\}
\newtheorem{Theorem}[Lemma]{Theorem}\newcommand{\bthe}{\begin{Theorem}}\newcommand{\ethe}{\end{Theorem}}
\newtheorem{Definition}[Lemma]{Definition}\newcommand{\bed}{\begin{Definition}}\newcommand{\eed}{\end{Definition}}
\newtheorem{Remark}[Lemma]{Remark}\newcommand{\beR}{\begin{Remark}\rm}\newcommand{\eeR}{\end{Remark}}
\newtheorem{Corollary}[Lemma]{Corollary}\newcommand{\bcor}{\begin{Corollary}}\newcommand{\ecor}{\end{Corollary}}

\newcommand{\bet}{\begin{tabular}{lll}}\newcommand{\eet}{\end{tabular}}

\title[]{E\MakeLowercase{quivalence of germs (of mappings and sets) over $\k$ vs that over $\K$}}
\author[]{D\MakeLowercase{mitry} K\MakeLowercase{erner}}

\address{\tiny Department of Mathematics, Ben-Gurion University of the Negev, P.O.B. 653, Be'er Sheva 84105, Israel. dmitry.kerner@gmail.com}

\date{\today\ \  filename: \jobname.tex}
\thanks{I was supported by the Israel Science Foundation, grants No.  1910/18 and 1405/22}


\begin{document}
\maketitle
\vspace{-1cm}
\begin{abstract}
Consider $\R$-analytic mapping-germs, $(\R^n,o)\stackrel{f,\tf}{\to}(\R^m,o).$ They can be equivalent (by coordinate changes) $\C$-analytically, but not $\R$-analytically.
  However, if the transformation of $\C$-equivalence is $Id $ modulo higher order terms, then it implies the $\R$-equivalence.

  On the other hand, starting from $\C$-analytic map-germs $(\C^n,o)\stackrel{f,\tf}{\to}(\C^m,o),$ and taking any field extension $\C\hookrightarrow\K,$ one has: if $f\stackrel{\K}{\sim}\tf$ then $f\stackrel{\C}{\sim}\tf.$

\medskip

  These (quite useful) properties seem to be not well known.
 We prove slightly stronger properties in a more general form:
\bei
\item for $\Maps,$ where $X,Y$ are (formal/analytic/$\k$-Nash) germs of spaces, with arbitrary singularities, over a base ring $\k;$
\item for the classical groups of (right/left-right/contact) equivalence  of Singularity Theory;
\item for faithfully-flat extensions of rings $\k\hookrightarrow\K.$ In particular, for arbitrary extensions of fields.
\eei

\medskip

The case ``$\k$ is a ring" is important for the study of deformations/unfoldings. E.g. it implies the statement for field extensions: if a family of $\k$-maps $\{f_t\}$ is $\K$-trivial    then it is also $\k$-trivial.

Similar statements for germs of spaces (``isomorphism over $\K$ vs isomorphism over $\k$") follow by the standard reduction ``Two maps are contact equivalent \iff their zero sets are ambient isomorphic".

\medskip

This study involves the contact equivalence of maps with singular targets, which  seems to be not well-established. We write down the relevant part of this theory.
\end{abstract}

\setcounter{secnumdepth}{6} \setcounter{tocdepth}{1}
\vspace{-0.2cm} \tableofcontents
\vspace{-0.8cm}

\section{Introduction}
\subsection{}Consider mapping-germs $\Mapk.$ Here $\k$ is an arbitrary field, $(\k^n,o)$ is the (algebraic, resp.  analytic, formal) germ of the affine space.
 See \S\ref{Sec.Background.Schemes.Maps} for notations and conventions.

 These mapping-germs are traditionally studied up to the right ($\cR$), left-right ($\cL\cR$ or $\cA$), and contact ($\cK$) equivalences. The case  of analytic maps,
  $\MapR$ and $\MapC,$  was intensively studied through  20'th century, see \cite{Mond-Nuno}.

 Many results of $\cR,\cK$-equivalences were extended (in numerous works) from $\R,\C$ to arbitrary fields,
 see Chapter 3 of the second edition of \cite{Gr.Lo.Sh}. The extension of $\cL\cR$-equivalence to   arbitrary fields  was done in \cite{Kerner.Group.Orbits}.
 Below $\cG$ is one of the groups $\cL,\cR,\cL\cR,$ $\cC,\cK,$ see \S\ref{Sec.G.}.

\medskip

The classical approach is: first to study the complex-analytic situation, then to deduce results in the real-analytic case.
 Hence the old question is to relate  the real and complex equivalence  of real map-germs.   Suppose  $f,\tf\sin \MapR$ are $\cG_\C$-equivalent, i.e. their complexifications are equivalent as elements of $\MapC.$ This does not imply their $\cG_\R$-equivalence.  The simplest example is the Morse critical point, $\sum(\pm)x^2_i$. Indeed, the real right equivalence, $\cR_\R,$ preserves the signature of this quadratic form, $(\la_+,\la_0,\la_-).$ (Thus e.g.   $\sum x^2_i\stackrel{\cR_\R}{\not\sim}-\sum x^2_i.$) The real contact equivalence, $\cK_\R,$ can swap the signs of all the eigenvalues,
  $(\la_+,\la_0,\la_-)\leftrightarrow(\la_-,\la_0,\la_+).$ But the difference $|\la_+-\la_- |$ is still preserved.
However one might guess:
\[\text{``The real and complex equivalence of real mapping-germs can differ only at low orders."}
\]
Namely, if $\tf=g_\C f,$ where the element $g_\C\sin \cG_\C$ is identity modulo   higher order terms, then $f\stackrel{\cG_\R}{\sim}\tf.$

\medskip

Even if one studies $\MapC,$   occasionally one has to extend the base-field. For example, some proofs go via deformation arguments, i.e. additional parameters are introduced, $\C\hookrightarrow \C(t)=:\K.$ Then, having established the $\cG_\K$-equivalence, one would like to deduce the $\cG_\C$-equivalence. The natural guess is: ``This should work, because $\C$ is algebraically closed".

\subsection{The results}\!We convert these guesses into stronger statements in  more general situation:
 \bei
 \item
 For $\Maps$   whose source $X\sseteq(\k^n,o)$ and target $Y\sseteq (\k^m,o)$ are (formal or analytic or algebraic) germs over a ring $\k,$ and with arbitrary singularities.
 \item
For a large class of ring extensions $\k\sset\K.$ (In particular for arbitrary field extensions.)
 \eei

The precise statement is in  Theorem \ref{Thm.Equivalence.Change.base.ring}, for faithfully flat ring extensions. Now we give just the conclusion for (arbitrary) field extensions $\k\sset \K.$
 \bee[\bf 1.]
 \item (For  $char(\k)=0$) \quad  If $\tf\stackrel{\cG_\K}{\sim}f$ by an element whose linear part is unipotent, then
   $\tf\stackrel{\cG}{\sim}f.$
 \item  Suppose the base-field is algebraically closed, of any characteristic.
 \mbox{If $\tf \stackrel{\cG_\K}{\sim} f $ then
   $\tf\!\stackrel{\cG}{\sim}\!f.$}
 \eee

Part 1 is sharp, as we show by explicit examples in \S\ref{Sec.Remarks.Sharpness}. While part 2 is ``naturally expected", part 1 seems ``a bit unexpected". (At least for some people in Singularities.)

\medskip

\noindent Below we list some corollaries of Theorem \ref{Thm.Equivalence.Change.base.ring} and other results.

\bei
\item
Working over base rings (not just base fields)  is important for deformations and unfoldings. E.g. it implies,
 Corollary \ref{Thm.Equivalence.for.families}, for $char(\k)=0$:
\\{\em Let $\k\sset \K$ be any extension of fields. If a family $f_t$ is $\cG_\K$-trivial, then it is also $\cG$-trivial.}

(Roughly: the $\k$-part of the $\cG_\K$-constant stratum locally coincides with the $\cG_\k$-constant stratum.)
 This seems to be new even in the classical case of $\R\sset\C.$
\\{\bf Example \ref{Ex.Stability.base.change}.} {\em $f\in \Maps$ is stable as a $\k$-map \iff it is stable as a $\K$-map.}

\medskip

\item The ``$\k$ vs $\K$ question" is finitely determined. Namely, if two maps $f,\tf\sin \Maps$ are $\cG_\K$-equivalent, then their $\cG$-equivalence can be verified at a finite jet.   (Here $f,\tf$ are not assumed finitely determined.)
    More precisely, Corollary \ref{Thm.K.vs.k.finitely.determined}  reads   in $char=0$:

\bee
\item  Every map $f\sin \Maps$ satisfies: $\cG_\K f\cap (\{f\}+\cm^N\cdot \RmX)\sset \cG f$ for some $N\sin \N.$
\item Moreover, this number $N$ can be chosen uniformly for the whole $\cG_\K$-orbit of $f$.
\eee
   \medskip

 The smallest such  $N$ is of course crucial for various applications, but the proof is based on Artin's $\be$-function, for which not many bounds are known.

\medskip

\item
Though the $\cG_\K$-equivalence does not imply the $\cG$-equivalence, we have at least the finiteness of $\cG$-orbits.
\\{\bf Corollary \ref{Thm.Finite.Splitting.orbits}.} {\em The real part of the complex orbit, $\cG_\C f\cap \Maps,$ splits into a \underline{finite} number of real orbits, $\amalg \cG_\R f_\al.$}
\\In fact Corollary \ref{Thm.Finite.Splitting.orbits} is for finite extensions of ``fields of type $F$'', e.g. $\R$ or $p$-adic fields.

\medskip

\item The classical group-actions on spaces of map-germs, $\cG\circlearrowright\MapX,$ for $\cG=\cL,\cR,$ $\cL\cR, \cC,\cK,$ are well-studied.
 Maps with singular targets, $\Maps,$ were studied only in rather particular cases, e.g. \cite{Mond-Montaldi}, \cite{Nun-Ball.Tomazella}.
 We could not find the general references for the actions
 $\cG\circlearrowright\Maps$  with singular targets.  The relevant basic notions are set-up in \S\ref{Sec.G.and.TG}.

Giving the full exposition of this topic    ``$\cG\circlearrowright\Maps$" would (at least) double the paper. Therefore only the main definitions are given, and precise references for all the ingredients are provided.

  We define the tangent spaces to the groups at unit element, $T_\cG.$ In full analogy with the classical ``Lie-algebra action on manifolds" one  gets the  action $T_\cG\times\Maps\to T(Maps).$ For each map we get the
   tangent image, $T_\cG f\sseteq T_{(Maps,f)}.$ Specifying a filtration on $\Maps$ defines the corresponding filtrations of the group, $\cG^\bullet<\cG,$ and of the tangent space,
  $T_{\cG^\bullet}\sset T_\cG.$

\medskip

\item Two auxiliary tools are obtained: the $T_\cG$-version  of the Artin-Rees lemma,  \S\ref{Sec.Artin-Rees.TG},  and the $\cG$-version, \S\ref{Sec.Artin-Rees.G}. The latter seems to be the first non-linear version of  the Artin-Rees lemma in Commutative Algebra.
 Its immediate application is:
\\{\bf Example \ref{Ex.G.-fin.det.iff.Gjj.fini.det}.} {\em For any $f\in \Maps$ the following conditions are equivalent:
\bee
\item $f$is $\cG$-finitely determined
\item  $f$ is $\cG^{(j)}$-finitely determined for some $j\sin \N $
 \item  $f$ is $\cG^{(j)}$-finitely determined for every $j\sin \N.$
\eee
}

This is well known in $char=0,$ where such problems are transformed to the tangent level by the exponential/logarithmic maps.
 In $char>0$ this statement seems to be new.

\medskip

\item While most of our results are local (for map-germs), in \S\ref{Sec.K.vs.k.Global} we give an immediate global implication,
 in $char=0$. {\em If a family of projective morphisms $\P X\stackrel{\P f_t}{\to }\P Y $   is globally-$\K$-trivial, then it is globally-$\k$-trivial.}
\eei

\subsection{Applications}
Besides the applications in \S\ref{Sec.K.vs.k.Examples.Corollaries}, these results are of immediate use, e.g. in real Singularity Theory and real Algebraic/Analytic Geometry. We recall just two particular directions:
\bei
\item (Hilbert 17'th problem) representing real polynomials or analytic/meromorphic functions as sums of squares. See e.g. \cite{Fernando.24}, \cite{Blekherman.Smith.Velasco} and further references.
\item (Nash-Tognoli problems) Approximations of $\R$-analytic sets/functions by their Nash counterparts. See e.g. \cite{Ghiloni-Savi}, \cite{Bochnak.Coste.Roy}
 and further references.
\eei
In both cases an essential task is to pass from $\R$ to $\C,$ and back to $\R.$

\

 Our results seem helpful in addressing the old and open question ``When/whether the real-$\cR$-modality of function germs coincide with their complex-$\cR$-modality?", \cite[pg.34]{AGLV}.

 Of course, the similar $\cL\cR$-question for maps is even more interesting.

 Another particular motivation comes from the results of  Mather-Yau/Gaffney-Hauser-type, \cite{Kerner.Unfoldings}. In those proofs one has to pass to the algebraic closure, $\k\sset \bk,$ but at the end we should deduce the equivalence over the initial field $\k.$

\subsection{Remarks}
\bee[\bf i.]
\item Our results go in the style of the classical Mather's lemma, \cite[Theorem 6.3]{Mond-Nuno}, and of the ``unipotence ideology" of \cite{Bruce.duPlessis.Wall}.

Even for the case of ``$\MapR$ vs $\MapC$" the  results seem to be  not well known.
   A particular case of Theorem \ref{Thm.Equivalence.Change.base.ring} has appeared in the first version of \cite{Kerner.Unfoldings}.

\item Suppose the target is   $\k$-smooth, i.e. $Y\=(\k^m,o).$  Then the contact equivalence of map-germs translates into the  ambient isomorphism of their zero sets, more precisely, of the defining ideals.
    Thus the geometric interpretation of the $\cK$-part of Theorem \ref{Thm.Equivalence.Change.base.ring} is:
     ``If two germs of spaces are unipotent-isomorphic over $\K,$ then they are unipotent-isomorphic over $\k$".

    The left-right equivalence of map-germs translates into   isomorphism  of the corresponding foliation-germs.
        Thus the geometric interpretation of the $\cL\cR$-part of Theorem \ref{Thm.Equivalence.Change.base.ring} is:
   ``If two germs of foliations are unipotent-isomorphic over $\K,$ then they are unipotent-isomorphic over $\k$".

  \item Recall the old ``inverse question": Suppose $f,\tf\sin \MapC$ are $\R$-analytically equivalent (as maps $(\R^{2n},o)\to(\R^{2m},o)$),
   does this imply the $\C$-analytic equivalence? See \cite{Ephraim} and \cite{Fern.-Hern.Gim.-Con.}.
    In the subsequent paper we apply our results to address this question.
\eee

 \subsection{Acknowledgements} Part of this work was done during  the ``Thematic Month on Singularities",
   Luminy, January-February 2025.
\\
The results of the papers \cite{Kosar.Popescu}, \cite{Popescu.Rond} were especially useful.
\\
Thanks are to Mikhail Borovoi, for the help on Galois cohomology and on ``fields of type $F$". (In particular for Corollary \ref{Thm.Finite.Splitting.orbits} and Remark \ref{Rem.Finite.Splitting.of.Orbits}.ii.)

We thank referee for important remarks.

\section{Notations  and conventions}\label{Sec.Background.Schemes.Maps}
 For the general introduction to  $\Mapk$
 in the real/complex-analytic case see \cite{Mond-Nuno}.  For some definitions and results on  $\Maps$ over arbitrary fields
  see \cite{Kerner.Group.Orbits}.

\medskip

  We use   multivariables, $x=(x_1,\dots,x_n),$  $y=(y_1,\dots,y_m),$ $t=(t_1,\dots,t_l),$ $F=(F_1,\dots,F_r).$
 The notation  $F(x,y)\sin (x,y)^N$  means: $F_j(x,y)\sin (x,y)^N$ for each $j.$

 \subsection{The germs of spaces}\label{Sec.Background.Germs.Rings} The source  and the target of our maps, $X,$ $Y,$ are (formal, resp. analytic,  algebraic) germs of spaces. Here:
\bei
\item  $X\!=\!\Spec(R_X)$, for the rings  $R_X\!=\!\quots{\k[\![x]\!]}{J_X},$ resp. $R_X\!=\!\quots{\k\{x\}}{J_X},$   $R_X\!=\!\quots{\k\bl x\br}{J_X}.$
  \\$Y=\Spec(R_Y)$, for the rings $R_Y\!=\!\quots{\k[\![y]\!]}{J_Y},$ resp. $R_Y\!=\!\quots{\k\{y\}}{J_Y},$   $R_Y\!=\!\quots{\k\bl y\br}{J_Y}.$
  \\Accordingly $R_{X\times Y}$ is one of   the rings $ \quots{\k[\![x,y]\!]}{(J_X+J_Y)},$ $\quots{\k\{x,y\}}{(J_X+J_Y)},$   $ \quots{\k\bl x,y\br}{(J_X+J_Y)}.$

\item For formal power series, $R_X=\quots{\k[\![x]\!]}{J_X},$ the base ring $\k$ is either a field, or $\k=\quots{\k_o[\![t]\!]}{\ca},$
 for an arbitrary field $\k_o.$

For analytic power series, $\quots{\k\{x\}}{J_X},$   either  $\k$ is a normed field   or
 $\k=\quots{\k_o\{t\}}{\ca},$ where $\k_o$ is  a normed field.
 The norm is always non-trivial, and the field is always assumed quasi-complete \wrt its norm\footnote{A normed field, $\k,\|*\|,$ is called quasi-complete, if the $\|*\|$-completion,   $\bk^{\|*\|},$ is separable over $\k.$ Thus any normed field of zero characteristic  is quasi-complete, and so is any complete field.  The norm is always assumed non-trivial.}. The simplest examples are $\Q,\R,\C$ or the $p$-adic field $\Q_p.$


    For algebraic power series, $\quots{\k\bl x\br}{J_X},$   $\k$ is any  field or an excellent Henselian local ring over a field.
  The main relevant (non-field) example is $\k=\quots{\k_o\bl t\br}{\ca},$ for $\k_o$ a field.
 See page 4 of \cite{Denef-Lipshitz} for more detail.

 The ring $\k\bl x\br$ consists of power series satisfying polynomial equations over $\k[x].$ E.g. for $char(\k)=0$ the Taylor expansion of $(1+x)^\al$ (with $\al\sin \Q$) is an algebraic power series. But the power series of $e^x,$ $ln(1+x)$ are non-algebraic.
  The ring $\R\bl x\br$ is just the ring of germs of Nash functions.
 Therefore we call the case  $\quots{\k\bl x\br}{J_X},$ as ``$\k$-Nash" (or just ``Nash"), to avoid any confusion with the algebraic rings like  $\k[x]$ or  $\k[x]_{(x)}.$

\medskip

In all these cases we call $\k$ (resp. $\k_o$) {\em the base field}.

 \item

  In the $\k$-smooth case (i.e.  $J_X=0,J_Y=0$) we denote  $(\k^n,o):=\Spec(R_X)$, $(\k^m,o):=\Spec(R_Y).$ These are formal/analytic/$\k$-Nash germs of spaces,      depending on the context.

In the general case we have the subgerms $X=V(J_X)\sseteq (\k^n,o)$ and $Y=V(J_Y)\sseteq (\k^m,o),$ with the structure rings $R_X,R_Y.$
The singularities of the germs $X,Y$ are arbitrary, not necessarily isolated, or complete intersections, or reduced.
\item
   Denote the maximal ideal of a ring $\k$  by $\cm_\k.$ (If $\k$ is a field, then $\cm_\k=0$.)

Fix some coordinates $(x_1,\dots,x_n)$ on $(\k^n,o),$ and $(y_1,\dots,y_m)$ on  $(\k^m,o).$
 They are sent to the generators of ideals in $R_X,$ $R_Y.$  Abusing notations we write just  $(x)=(x_1,\dots,x_n)\sset R_X$ and $(y)=(y_1,\dots,y_m)\sset R_Y.$
     If $\k$ is a field, then $(x),(y)$ are the maximal ideals. Otherwise the maximal ideals are $\cm_X:=\cm_\k+(x)\sset R_X$
      and $\cm_Y:=\cm_\k+(y)\sset R_Y.$  When only $R_X$ is involved we write just $\cm\sset R_X.$

      \item
 The rings $R:=\k[\![x]\!],\k\{x\},\k\bl x\br$ admit substitutions: if $f\sin R$ and $g_1,\dots,g_n\sin (x)\sset R,$ then $f(g_1,\dots,g_n)\sin R.$  This fact is trivial in the formal case, and well-known in the analytic case, \cite[pg.70]{Serre.Lie}. For $\k$-Nash case see Remark 1.3 of \cite{Denef-Lipshitz}.
\eei

\beR\label{Rem.Noetherian.Submodules}
All our rings are Noetherian.
  Recall the general fact: any $R_X$-submodule $M\sseteq \RmX$ is finitely generated. Proof: $\RmX$ itself is a f.g. module over the Noetherian ring $R_X$,
 therefore $\RmX$ is a Noetherian module. Hence $M$ is f.g.
\eeR

\subsection{The module of derivations}\label{Sec.Background.Derivations}
 The study of germs of sets/mappings traditionally goes via vector fields. Algebraically these are $\k$-linear derivations of the structure ring,
  $Der_X:=Der_\k R_X.$ See e.g. \cite[\S2.2]{Mond-Nuno}, \cite[\S1.1.10]{Gr.Lo.Sh}, \cite[\S25]{Matsumura}.

 \bex\label{Ex.Derivations}
 \bee[\bf i.]
 \item For $X=(\k^n,o)$ one gets the free $R_{(\k^n,o)}$-module of rank=$n$,  i.e. $Der_{(\k^n,o)} =R_{(\k^n,o)}\bl\frac{\di}{\di x_1},\dots,\frac{\di}{\di x_n}\br.$
 See e.g. Theorem 30.6(ii) of \cite{Matsumura}.
 \item
 In the general case, i.e.
   $X=V(J_X)\sset(\k^n,o),$ the elements of $Der_X$  lift to logarithmic derivations on $(\k^n,o)$ that preserve the ideal $J_X.$
 Namely, the following sequence is exact:
\beq
0\to J_X\cdot Der_{(\k^n,o)}\to Der_{(\k^n,o)}Log(J_X)\to Der_X\to0,
\eeq
   where $Der_{(\k^n,o)}Log(J_X):=\{\xi\sin Der_{(\k^n,o)}|\ \xi(J_X)\sseteq J_X\}.$
 \eee
 \eex

\noindent The $R_X$-module $Der_X$ is finitely-generated. Indeed, it is the quotient of the module $Der_{(\k^n,o)}Log(J_X),$ the latter being finitely generated by
 Remark \ref{Rem.Noetherian.Submodules}.

Applying the module of derivations to an element $f\in R_X$ we get the ideal $Der_X(f)\sseteq R_X.$ For an element $f\in \RmX$ we get the submodule
 $Der_X(f)\sseteq\RmX.$

   \subsection{Maps of scheme-germs}
A (formal/analytic/Nash) map $f\sin \Maps$ is defined by the corresponding homomorphism of $\k$-algebras, $f^\sharp\sin Hom_\k(R_Y,R_X).$ All our homomorphisms are local.

\bei
\item In the   case   $J_Y=0 $    one has  maps to a $\k$-smooth target,  $\MapX.$ We write this explicitly.
  With the chosen generators $y=(y_1,\dots,y_m)\sset R_Y$  each map $X\stackrel{f}{\to} (\k^m,o)$ is presented by
 a tuple of (formal/analytic/algebraic) power series, $f:=(f_1,\dots,f_m)\sin \cm_X\cdot\RmX.$
 Vice versa, every element of $\cm_X\cdot\RmX$ defines a map $X\stackrel{f}{\to}Y.$ This identifies the space of maps,
  $\MapX=\cm_X\cdot \RmX.$

\item
For an arbitrary target   $Y\sseteq (\k^m,o),$   with $R_Y=\quots{\k[\![y]\!]}{J_Y},\quots{\k\{y\}}{J_Y},\quots{\k\bl y\br}{J_Y},$   a map $X\stackrel{f}{\to} Y$ is still represented by
  an element $f \sin \cm_X\cdot\RmX.$ But now the condition $f^\sharp\sin Hom_\k(R_Y,R_X)$ means: $f^\sharp(J_Y)\sseteq J_X\sset R_{(\k^n,o)}.$
  This embeds the space of maps as a subset:
 \beq
 \Maps=Hom_\k(R_Y,R_X)=\{f\sin \cm_X\cdot \RmX\ |\quad f^\sharp(J_Y)\sseteq J_X\}.
 \eeq

\noindent
For  not $\k$-smooth  targets     the subset $\Maps\sset \cm_X\scdot \RmX$ is not closed under addition.
\eei

\subsection{Filtrations on $\Maps$}\label{Sec.Background.Filtrations.of.Maps}
Take a filtration by $\k$-submodules \mbox{$\RmX\!=\!\tM^0\spset\tM^1\spset\cdots.$}
 Suppose this filtration is equivalent to the filtration $I^\bullet\cdot \RmX,$ for some ideal $I\sset  R_X.$
  Namely:
  \beq
  \tM^j\spset I^{d_j}\cdot \RmX \quad \text{ and } \quad I^j\cdot \RmX\spset \tM^{d_j}, \quad \text{ for each } j \text{ and a corresponding } d_j\gg j.
  \eeq

Take the corresponding filtration on the space of maps, $M^\bullet:=\Maps\cap \tM^\bullet.$
 Explicitly: $M^j=\{f\sin \tM^j|\ f^\sharp(J_Y)\sseteq J_X\}.$
 \bex
 \bee[\!\!\bf i.]
 \item
  For $\k$-smooth targets, $Y\!=\!(\k^m,o),$ we get just the submodules    $M^j\!=\!  \tM^j\sset \RmX.$

\item For any $Y\sseteq(\k^m,o)$ and $\tM^\bullet=\cm^\bullet\scdot\RmX$ one has:

 $M^{j}=\Maps\cap \cm^j\scdot\RmX=\{f\in Hom(R_Y,R_X)\ |\quad f(\cm_Y)\sseteq\cm^j_X \}.$
\eee
\eex
 The $M^\bullet$-order of maps is defined in the standard way, $ord(f):=sup\{d|\ f\sin M^d\}.$

\medskip

In many applications one needs just the filtration  $I^\bullet\cdot \RmX.$ However we work with the slightly more general version, $\tM^\bullet,$ as it is used in the proof of Corollary \ref{Thm.Equivalence.change.base.field}.
 In particular, we do not assume $\tM^1\sseteq \cm\cdot \RmX.$

\subsection{The tangent space to $\Maps$}\label{Sec.Tangent.to.Maps} In \S\ref{Sec.TG.action.on.Maps} we work with tangent spaces of various groups of equivalence. Hence we need the tangent space  to $\Maps$ at an element $f,$ denote it $T_{(Maps,f)}.$ It is obtained in the standard way, \cite[pg.80]{Hartshorne}.
 Fix a set of generators $(q_1,\dots,q_r)$ of $J_Y\sset R_{(\k^m,o)},$ these are power series.  Take the  matrix of their derivatives $[q'_\bullet]\!\in\! Mat_{r\times m}(R_{(\k^m,o)}).$ Evaluate its entries at $f,$ we get the matrix $[q'_\bullet|_{y=f}]\!\in\! Mat_{r\times m}(R_X).$
   Then
   \beq\label{Eq.Tangent.to.Maps.presentation}
   T_{(Maps,f)}:=\{v|\quad [q'_\bullet|_{f}]\cdot v=0\in \RrX\}=ker\big[\RmX\stackrel{[q'_\bullet|_{f}]}{\to}\RrX\big]\sseteq \RmX.
\eeq
 By the direct check, $T_{(Maps,f)}\sseteq \RmX$ is an $R_X$-submodule, and it does not depend on the choice of generators $q_\bullet.$
  By Remark \ref{Rem.Noetherian.Submodules} this submodule is finitely generated.
  \bex
  \bee[\bf i.]
\item For $J_Y=0,$ i.e. $Y=(\k^m,o),$ we get: $T_{(Maps,f)}= \RmX,$ a free module of rank $m.$ In particular, this tangent space is independent of $f.$
\item
For any ideal $0\neq J_Y\sseteq (y)^2,$ and the constant map, $f=0,$ we get $T_{(Maps,0)}= \RmX.$
 For non-constant maps we get in general $T_{(Maps,f)}\ssetneq \RmX.$
\eee
\eex

\subsection{Faithfully flat extensions, see e.g. \cite[Appendix B]{Gr.Lo.Sh}}\label{Sec.Background.FF.extensions} An extension of rings, $\k\sset \K,$ is called faithfully-flat (f.f.) if the following conditions hold  for any  (not necessarily finitely-generated) $\k$-modules:
\bei
\item (faithfulness) If $V\neq0$ then $\K\otimes V\neq0.$
\item (flatness) If the sequence $0\to V\to U\to W\to 0$ is exact, then the sequence   $0\to \K\otimes V\to \K\otimes U\to\K\otimes W\to 0$ is exact.
\eei
\bex\label{Ex.Faithfully.Flat}
\bee[\bf i.]
\item For $\k$ a field, any ring extension $\k\sset \K$ is f.f.
\item For any ring $\k$ the extension $\k \sset \k[\![ x]\!]$ is f.f.
 And similarly for  $\k \sset \k\{ x\}$ and $\k \sset \k\bl  x\br.$

\item For any field extension $\k_o\sset \K_o$ the extension of rings $\quots{\k_o[\![t]\!]}{\ca}\sset \quots{\K_o[\![t]\!]}{\ca}$ is f.f.
\eee
\eex
 Take a f.f. extension of rings $\k\sset \K.$ Then the natural map $V\to \K\otimes V,$ by $v\to 1\otimes v,$ is an embedding.
  Below we identify $V$ with its image $1\otimes V\sset \K\otimes V.$ A $\k$-submodule $V\sset W$ leads to a $\K$-submodule $\K\otimes V\sset \K\otimes W.$ We get the obvious inclusion  $(\K\otimes V)\cap W\supseteq V.$

 \bel\label{Thm.Faithful.Flatness}
If the extension  $\k\sset\K$ is faithfully-flat, then $( \K\! \otimes \!V )\!\cap\! W\! =\!V $.
 \eel
This is used in the proof of Theorem \ref{Thm.Equivalence.Change.base.ring}.
 \bpr
 Denote $W':=( \K \otimes V )\cap W  ,$ and consider the exact sequence
  $0\to V\to W'\to \quots{W'}{V}\to0.$ Applying $\K\otimes$ we get   the exact sequence
   $0\to   \K\otimes V\to \K\otimes  W'\to \K\otimes \quots{W'}{V}\to0,$
     by flatness.
   Here the map    $\K\otimes V\to  \K\otimes  W'$ is an isomorphism.  Therefore
   $\K\otimes \quots{W'}{V}=0.$ The faithfulness of $\K\otimes$ gives then: $\quots{W'}{V}=0.$
\epr

 \section{The groups $\cR,\cL,\cL\cR,\cC,\cK,$ their filtrations, and their (filtered) tangent spaces}\label{Sec.G.and.TG}

For the case of $\Mapk$ the  groups   $\cR,\cL,\cL\cR,\cC,\cK,$  their tangent spaces $T_ \cR ,T_ \cL ,T_{ \cL\cR },T_\cC,T_{\cK },$ and their actions on  $\Mapk$  are classical.  See e.g. \S3.2 and \S4.4  in \cite{Mond-Nuno}
  in the case $\k=\R,\C$, and Definition 3.1.1. in \cite{Gr.Lo.Sh} for  the groups $\cR,\cK$ and an arbitrary field $\k$.
  For the filtered case, i.e. the groups $\cR^{(j)},\cL^{(j)},(\cL\cR)^{(j)},$ $\cC^{(j)},\cK^{(j)}$, over arbitrary $\k,$ acting on $\MapX,$
     see \cite[\S3]{Kerner.Group.Orbits}.

      The case of singular target, $Y\ssetneq(\k^m,o),$ seems to be not well established even for (real/complex)  analytic maps. We introduce the basic notions.   In the cases $J_X=0,J_Y=0$ this gives the classical objects.

\noindent
Below $R_X$ is one of the rings $ \quots{\k[\![x]\!]}{J_X} , \quots{\k\{x\}}{J_X} ,  \quots{\k\bl x\br}{J_X} .$
 Accordingly   $R_Y$ is one of $\quots{\k[\![y]\!]}{J_Y} ,$   $\quots{\k\{y\}}{J_Y} ,$   $\quots{\k\bl y\br}{J_Y} .$

\subsection{The groups $\cR,\cL\cR,\cC,\cK$}\label{Sec.G.}
 The  (formal/$\k$-analytic/$\k$-Nash) automorphisms of $X$ over the base $\Spec(\k)$ are defined via their algebraic counterparts, $\k$-linear automorphisms of the structure ring, $Aut_X:=Aut_\k(R_X).$
 Similarly we denote  $Aut_Y:=Aut_\k(R_Y).$
\bex
\bee[\bf i.]
\item In the $\k$-smooth case, $\Mapk,$ these automorphisms are the coordinate changes in the source and the target.
\item In the general case any automorphism $\Phi_X\in Aut_X$ is presentable as a coordinate change of the ambient germ $(\k^n,o)$ that preserves the source $X\sset (\k^n,o).$ And similarly for $\Phi_Y\in Aut_Y.$
\eee
\eex

\medskip

\noindent These automorphisms act on  maps, defining the standard
 group actions of Singularity Theory:
\bei
\item (the right group) $\cR\!:=Aut_{X}\!: =Aut_\k(R_X).$ The action $\cR\!\circlearrowright\! \Maps$ goes by $(\Phi_X,f)\!\to\! f\circ \Phi^{-1}_X.$
\item (the left group)  $\cL\!:=Aut_{Y}\!: =Aut_\k(R_Y).$ The action $\cL \!\circlearrowright\! \Maps$ goes by $ (\Phi_Y,f)\!\to\! \Phi_Y\circ f$.
\item (the left-right  or $\cA$-group) The action $  \cL \cR\!:=\cL\times\cR \!\circlearrowright\! \Maps $  goes by $ (\Phi_Y,\Phi_X,f)\to \Phi_Y\circ f\circ \Phi^{-1}_X.$
\item (the $y$-part of the contact group) The subgroup $\cC< Aut_{X\times Y}$ consists of automorphisms of $X\times Y$ that preserve the fibres of the projection $X\times Y\to X$
 and preserve ``the zero-level slice" $X\times\{o\}\sset X\times Y.$ Therefore $\cC$ acts as $Id_X$ on $X\times\{o\},$
  and restricts to automorphisms of the central fibre $\{o\}\times Y.$
\\
Explicitly, fix some generators $J_Y\!=\!\bl q_\bullet(y)\br\sset R_{(\k^m,o)},$ then the action $\cC\!\circlearrowright \!X\!\times\! Y$ is presentable by $(x,y)\!\to\! (x,C(x,y)),$
   where $q_\bullet(C(x,y))\sin R_{X\times Y}\cdot J_Y,$ and $C(x,y)\=\Phi_Y(y)+(\dots),$ with  $\Phi_Y\sin Aut_Y$ and  $(\dots)\sin \cm_X\scdot \cm_Y.$
\\The action of this group on maps, $\cC\!\circlearrowright\!\Maps,$ is defined via the graph, by $f(x)\!\rightsquigarrow \!C(x,f(x)).$
\item The contact group is the semi-direct product $\cK:=\cC\rtimes Aut_X.$  Note that $\cC\lhd\cK,$ a normal subgroup.
 The action $\cK\circlearrowright X\times Y$ goes by $(x,y)\to (\Phi_X(x),C(x,y)).$
 The action $\cK\circlearrowright\Maps$ is defined via the graph of $f,$ i.e. $f\rightsquigarrow C(x,f\circ \Phi^{-1}_X).$

\item
When the target is $\k$-smooth,  $Y=(\k^m,o),$ the contact group-action is ``reduced" to its linearized version, $\cK_{lin}=GL_m(R_X)\rtimes Aut_\k(R_X).$ Namely, the orbits of the two groups coincide: $\cK f=\cK_{lin}f.$
     See \cite[Proposition 4.2]{Mond-Nuno} for $f\sin \Mapk$ with $\k=\R,\C,$ and \cite[\S3.1]{Kerner.Group.Orbits} for   $f\sin \MapX.$

 When the target is singular, $Y\ssetneq(\k^m,o),$ the reduction of $\cK$ to $\cK_{lin}$  is more subtle, see \cite[\S2]{Kerner.G-AP}.
\eei

\medskip

Below $\cG$ is one of these groups, $\cL,\cR,\cL\cR,\cC,\cK,$ and we get the  traditional  notions of (left, right, left-right, contact) equivalence of maps, $ f\stackrel{\cG}{\sim }\tf.$

\subsection{Filtrations on the groups, $\cG^{(\bullet)}\le \cG$}\label{Sec.G.filtered}
 Take the filtration   $\tM^\bullet\sset \RmX$ of   \S\ref{Sec.Background.Filtrations.of.Maps}, and the corresponding filtration $M^\bullet\sset \Maps.$  It defines the subgroups of the groups $\cL,\cR,\cC$ for $j\!\ge\!0,$
 \beq\label{Def.filtered.groups}
 \cG^{(j)}:=\{g\sin \cG|\ g\cdot M^d=M^d  \text{ and }Id=[g]\circlearrowright \quot{M^d+\tM^{d+j} }{\tM^{d+j} } \quad \forall\ d\ge1 \}.
 \eeq
Thus we get the filtration $\cG\ge \cG^{(0)}\ge \cG^{(1)}\ge\cdots.$ See \cite[\S3.3]{Kerner.Group.Orbits} and \cite[\S2.2]{BGK.20} for   detail.

\bex\label{Ex.Filtered.Groups}
 Suppose  the target is $\k$-smooth, $Y=(\k^m,o).$ Thus $\Maps=\cm_X\cdot\RmX.$
 Take the filtration $M^\bullet=\cm^\bullet_X\cdot \RmX$ and some $j\ge1.$   We have, see \cite[\S3.3]{Kerner.Group.Orbits}:
\bei
\item $\cR^{(j)} \!=\!\{\Phi_X\sin Aut_X|\ \Phi_X(x)\!-\!x\sin \cm^{j+1}_X\},$ i.e. automorphisms of $X$ that are $Id$ modulo $\cm^{j+1}_X.$

\item For  $\k\spseteq\Q:$ \  $\cL^{(j)} =\{\Phi_Y\sin Aut_Y|\ \Phi_Y(y)-y\sin \cm^{j+1}_Y\},$ i.e. automorphisms of $Y$ that are $Id$ modulo $\cm^{j+1}_Y.$  See example 3.2 in \cite{Kerner.Group.Orbits}.

Similarly  $\cC^{(j)} =\{(Id_X,C)|\ C(x,y)-y\sin \cm_Y\cdot (\cm_X+\cm_X)^j \}.$
\item
  The group $\cC$ can be replaced by the linear group $GL_m(R_X),$  as in \S\ref{Sec.G.}. Then we have: $GL_m(R_X)^{(j)}=\{\one\}+Mat_{m\times m}(\cm^{j}_X).$
\eei
\eex

\noindent For the ``composite" groups  $\cL\cR,$ $\cK$ we \underline{define}:
\beq
(\cL\cR)^{(j)}:=\cL^{(j)}\times \cR^{(j)} \qquad  \text{ and } \qquad \cK^{(j)}:=\cC^{(j)}\rtimes \cR^{(j)}.
\eeq
\beR\label{Rem.on.FIltrations.of.G}
\bee[\bf i.]
\item  If one defines the filtration on the groups $\cL\cR,\cK$ directly via \eqref{Def.filtered.groups}, then the so obtained groups could be bigger than $\cL^{(j)}\times \cR^{(j)}, \cC^{(j)}\rtimes \cR^{(j)}$, e.g. when $\k$ is a finite field. But for $\k$  an infinite field, and $I^\bullet\neq 0$ for  each $\bullet,$ these versions coincide.
 See \cite[Lemma 3.4]{Kerner.Group.Orbits}. We remark that in Theorem \ref{Thm.Equivalence.Change.base.ring}
   the ring $\k$ is   infinite, being  an algebra either over  a field of $char=0$ or over an algebraically closed field.

\item

For the filtration $\tM^\bullet=\cm^\bullet_X\cdot \RmX$ one has  $\cG^{(0)}=\cG,$ as all the automorphisms are local.
For $\tM^\bullet=I^\bullet_X\cdot \RmX$ with $\sqrt{I}\ssetneq \cm_X$ we have  $\cG^{(0)}\lneq\cG.$

\item
These subgroups are normal, $\cG^{(j)}\lhd \cG^{(0)},$ see \cite[\S2.2.1]{BGK.20}.
\eee
\eeR

\subsection{The (extended) tangent spaces to the groups at unit elements, $\pmb{T_\cG}$}\label{Sec.TG.} These are defined in the standard way via  modules of derivations, as in \cite[\S3.2 and \S4.1]{Mond-Nuno}:
\beq
T_\cR:=Der_X,\hspace{2cm}T_\cL:=Der_Y,\hspace{2cm}T_{\cL\cR}:=T_\cL\oplus T_\cR.
\eeq

\bex
\bee[\bf i.]
\item Suppose the target is $\k$-smooth, i.e. $Y=(\k^m,o).$ Then $T_\cL$ is a free $R_Y$-module, $Der_Y=Der_{(\k^m,o)}=R_Y\bl\frac{\di}{\di y_1},\dots, \frac{\di}{\di y_m}\br.$
\item In the singular case, $X\!\ssetneq \!(\k^n,o),$ the tangent space $T_\cR$ is expressed via logarithmic derivations,  \S\ref{Sec.Background.Derivations}.
  \eee
\eex

\noindent
For the group $\cC<Aut_{X\times Y}$ we define $ T_\cC:=R_{X\times Y}\otimes_{R_Y} Der_Y Log(\cm_Y) \sset Der_{X\times Y},$ see \S\ref{Sec.Background.Derivations}
   The tensor product over $R_Y$ is well behaved as the module $Der_Y Log(\cm_Y)$ is finitely-generated.
\bex
\bee[\bf i.]
\item For $Y=(\k^m,o)$ the tangent space is: $T_\cC=R_{X\times Y}\scdot\cm_Y\scdot\bl\frac{\di}{\di y_1},\dots, \frac{\di}{\di y_m}\br=  \{\sum c_i(x,y)\frac{\di}{\di y_i}\},$
 for $c_\bullet(x,y)\in \cm_Y\scdot R_{X\times Y}.$
\item
When $J_Y\!\neq\!0,$ i.e. $Y\!\ssetneq\!(\k^m,o),$ we get:
$ T_\cC\=\{\sum \xi_i(x,y)\frac{\di}{\di y_i}\sin R_{X\times Y}\scdot Der_Y| \  \xi(\cm_Y)\sset R_{X\times Y}\scdot\cm_Y  \}. $
 \eee
\eex

Finally, the tangent space to the contact group is defined as $T_\cK:= T_\cC\oplus T_\cR.$

For $Y=(\k^m,o)$ we have also the linearized version,  $T_{\cK_{lin}}:=Mat_{m\times m}(R_X)\oplus T_\cR.$

\medskip

We have {\em defined} the tangent spaces $T_\cG.$ In \S\ref{Sec.TG.to.G.exp} this definition is justified  by the standard exponential/logarithmic maps
 $T_\cG\rightleftarrows\cG$ in $char=0.$ In \S\ref{Sec.TG.action.on.Maps} we give another link, the action $\cG\circlearrowright\Maps$ induces the action
  $T_\cG\times\Maps\to T(Maps).$

\subsection{The action  $\pmb{T_\cG\times\Maps\to T(\rm{Maps})}$}\label{Sec.TG.action.on.Maps}
Suppose a classical Lie group acts on a manifold, $G\circlearrowright Z.$
  This defines the Lie algebra action $T_{(G,1)}\times Z\to TZ,$ \cite[pg. 526]{Lee}. (Here $T_{(G,1)}$ is the tangent space at the unit element, while $TZ$ is the tangent bundle.) Indeed, for each vector $\xi\in T_{(G,1)}$ we have the flow $exp(t\xi).$ For each point $z_o\in Z$ take the orbit $\{exp(t\xi)z_o\}_t\sset Z.$
  Then the tangent vector to this orbit at $z_o$ is $\frac{d}{dt}(exp(t\xi)z_o)|_{t=0}=\xi z_o.$ In particular, $\xi z_o\in T_{(Z,z_o)}.$

\medskip

In the same way, for our groups $\cG$ the  elements of the tangent spaces define morphisms $T_\cG\times \Maps\to \RmX.$
 (In Lemma \ref{Thm.image.of.tangent.map} we show that in fact the image lies inside $T(Maps)$.)
 \bei
 \item
 $\pmb\cR\!:$ $(\xi_X,f)\to \xi_X(f):=(\xi(f_1),\dots,\xi(f_m))\in \RmX.$

 Here each $f_i$ is a power series ($mod\ J_X$) and $\xi(f_i)\in R_X$ is well defined, as $\xi(J_X)\sseteq J_X.$
 \item  $\pmb\cL\!:$ $(\xi_Y,f)\!\to\! \xi_Y(y)|_f\!:=\!(\xi_Y(y_1)|_f,\dots,\xi_Y(y_m)|_f)\sin \RmX.$
 \quad
 For the presentation $\xi\=\sum a_i(y)\frac{\di}{\di y_i}$ we get  $\xi(y_i)|_f\=a_i(f(x)).$ Here $\xi_Y$ is defined $mod\ J_Y\scdot Der_{(\k^m,o)},$ see Example \ref{Ex.Derivations}, while $f$ is defined $mod\ J_X\scdot\RmX.$ But the element $\xi_Y(y)|_f\sin \RmX$ is well defined, as $f^\sharp(J_Y)\sseteq J_X.$
\item
$\pmb\cC\!:$ $(\xi,f)\to \xi(y)|_f:=(\xi(y_1)|_f,\dots,\xi(y_m)|_f)\sin \RmX.$  For the presentation $\xi\=\sum a_i(x,y)\frac{\di}{\di y_i}$ we get  $\xi(y_i)|_f\=a_i(x,f(x))\sin R_X.$ The element $\xi(y)|_f\sin \RmX$ is well defined, as in the $\cL$-case.
\item For the composite groups, $\cL\cR,$ $\cK,$ the actions of $T_{\cL\cR}$ and $T_\cK$ are defined accordingly.
\eei

\bex
For $Y=(\k^m,o)$ we have $\MapX=\cm\cdot \RmX $ and $T_{Maps}=\RmX.$  One gets the tangent space actions   $T_\cG\times \cm\scdot \RmX\to \RmX.$
 See \cite[\S3]{Kerner.Group.Orbits} for more detail and examples. The case $\Mapk$ is classical, see e.g. \cite[\S3.2.1.]{Mond-Nuno}.
\eex

For a fixed map $X\stackrel{f}{\to}Y$ and a group $\cG$ we get the ``image tangent space" $T_\cG f\sseteq \RmX.$ Explicitly:
\bei
\item $T_\cR f\!=\!Der_X(f).$
\item $ T_\cL f\!=\!Der_Y(y)|_f\!=\!\{(\xi_{Y,1}|_f,\dots,\xi_{Y,m}|_f )|\ \xi_Y\!\in\! Der_Y\!\}.$
\item $T_\cC f\!=\!R_X\!\cdot\! Der_Y Log(\cm_Y)\cdot (y)|_f,$ see \S\ref{Sec.Background.Derivations}
\item
$T_\cK f\=T_\cR f+T_\cC f$ \ and \ $T_{\cL\cR}f=T_\cL f+T_\cR f.$
\item For $\k$-smooth target,  $Y=(\k^m,o),$  we have $ T_{\cK_{lin}}f=(f)\cdot \RmX+Der_X(f) .$
\eei
For the relation of these ``image tangent spaces", $T_\cG f,$ to the ``tangent spaces to the orbits", $T_{(\cG f,f)},$ see
 Theorem 3.1.21 of \cite{Gr.Lo.Sh}.

\medskip

The spaces $T_\cR f,T_\cC f,T_\cK f,T_{\cK_{lin}} f\sseteq \RmX$ are $R_X$-submodules. And they are finitely generated, by Remark \ref{Rem.Noetherian.Submodules}.
The spaces  $T_\cL f,T_{\cL\cR} f\sseteq \RmX$ are only $R_Y$-submodules.
\bex
Let $Y=(\k^1,o),$ then $T_{\cL \cR}f=Der_X(f)+f^\sharp(R_Y).$ Here $f^\sharp(R_Y)=\k\bl f\br$ (resp. $\k\{f\},$ $\k[\![f]\!]$) is the  ring of (algebraic, resp. analytic, formal) power series in $f.$
\eex

  Let $Y=(\k^m,o),$ then $T_{(Maps,f)}=\RmX.$ Thus $T_\cG f\sseteq T_{(Maps,f)}.$ This holds also in the singular case:

\bel\label{Thm.image.of.tangent.map}
Let $\cG=\cL,\cR,\cL\cR,\cC,\cK$ and $f\in \Maps.$ Then $T_\cG f\sseteq T_{(Maps,f)}\sseteq\RmX.$
 \eel
\bpr In view of the presentation of $T_{(Maps,f)},$ equation \eqref{Eq.Tangent.to.Maps.presentation}, it is enough to verify:
 $[q'_\bullet|_f]\cdot \xi(f)=0$ for each $\xi\in T_\cG.$ This is immediate by the chain rule:
 \bei
 \item
$\pmb\cR\!:$ \quad  $[q'_\bullet|_f]\cdot \xi_X(f)=\xi_X [q_\bullet|_f]\in \xi(J_X)\cdot \RmX\sseteq J_X\cdot\RmX=0.$
 \item
$\pmb\cL\!:$ \quad   $[q'_\bullet|_f]\cdot \xi_Y(y)|_f=[\xi_Y (q_\bullet)|_f]\in \xi_Y(J_Y)|_f\cdot \RmX\sseteq f^\sharp(J_Y)\cdot\RmX=0.$
\item $\pmb\cC\!:$  \quad  $T_\cC\sseteq R_{X\times Y}\otimes_{R_Y}T_\cL.$ Therefore $T_\cC f\sseteq  T_{(Maps,f)}.$

\item Similarly  $T_{\cL\cR}f=T_\cL f+T_\cR f\sseteq  T_{(Maps,f)}$ \quad and \quad $T_\cK f=T_\cC f+T_\cR f\sseteq  T_{(Maps,f)}.$
\epr
\eei

\subsection{Filtrations on the tangent spaces, $T_{\cG^{(\bullet)}}\sseteq T_\cG$}
Fix a  filtration $\tM^\bullet\sset \RmX,$ see  \S\ref{Sec.Background.Filtrations.of.Maps}.
   Accordingly the tangent spaces $T_\cG$  for $\cG=\cL,\cR,\cC $ acquire the filtration:
  \beq\quad
T_{\cG^{(j)}}:=\{\xi\sin T_\cG|\quad \xi(\Maps\cap \tM^d)\sseteq \tM^{d+j} \quad \forall\ d\ge1\},
\hspace{1cm}  T_\cG\spseteq  T_{\cG^{(0)}}\spseteq  T_{\cG^{(1)}}\spseteq\cdots
\eeq
For the composite groups $\cL\cR,\cK$ we {\em define} $T_{(\cL\cR)^{(j)}}:=T_{\cL^{(j)}}\oplus T_{\cR^{(j)}}$ and $T_{\cK^{(j)}}:=T_{\cC^{(j)}}\oplus  T_{\cR^{(j)}} .$
\bex  Fix an ideal $I \sset R_X$ and the corresponding filtration  $\tM^\bullet=I^\bullet\cdot \RmX.$
\bee[\bf i.]
\item
Let $Y=(\k^m,o),$ thus   $\Maps=\RmX.$ Therefore   $ M^\bullet=I^\bullet\cdot \RmX .$ Then the filtrations on the tangent spaces are (for $j\ge0$):
\bei
\item $T_{\cR^{(j)} }\supseteq I^{j+1}\cdot Der_X$ and $T_{\cL^{(j)} }\supseteq \cm^{j+1}_Y\cdot Der_Y.$
\item
$T_{\cK^{(j)}_{lin}}=Mat_{m\times m}(I^{j })\oplus T_{\cR^{(j)}}.$
\item
For $X=(\k^n,o)$ and $I=\cm_X$ we have: $T_{\cR^{(j)}}=\cm^{j+1}_X\cdot Der_{(\k^n,o)}.$
\item For $\k\spseteq\Q$  and $I=\cm_X$  we have: $T_{\cL^{(j)}}=\cm^{j+1}_Y\cdot Der_{(\k^m,o)}, $ see \S3.3.1 of \cite{Kerner.Group.Orbits}.

 Similarly we have  $T_{\cC^{(j)}}=(I+R_X\cdot \cm_Y)^j\cdot \cm_Y\cdot Der_{(\k^m,o)}.$
\eei
\item In the singular case, $Y\ssetneq (\k^m,o),$ the filtration  $T_{\cG^\bullet}\sset T_\cG$
 can be complicated. However we have immediate simple bounds.
 \bei
 \item $T_{\cR^{(j)}} \supseteq I\cdot T_{\cR^{(j-1)}} \supseteq\cdots \supseteq I^j\cdot T_{\cR^{(0)}}\supseteq I^{j+1}\cdot T_{\cR}.$
  \item $T_{\cL^{(j)}} \supseteq \cm_Y\cdot T_{\cL^{(j-1)}} \supseteq\cdots \supseteq \cm^j_Y\cdot T_{\cL^{(0)}}\supseteq \cm^{j+1}_Y\cdot T_{\cL}.$
   \item $T_{\cC^{(j)}} \supseteq (I+R_X\cdot \cm_Y)\cdot T_{\cC^{(j-1)}} \supseteq\cdots\supseteq   (I+R_X\cdot \cm_Y)^{j+1}\cdot T_\cC .$
 \eei
 \eee
\eex

\subsection{The exponential and logarithmic maps $T_\cG\rightleftarrows \cG$}\label{Sec.TG.to.G.exp}
Let  $R_X= \quots{\k[\![x]\!]}{J_X} $ and  $R_Y=\quots{\k[\![y]\!]}{J_Y},$ for  $\k\spseteq\Q.$  (Thus $char(\k)=0.$)
The relation between the groups and their tangent spaces is given   by the standard exponential/logarithmic maps. Namely:
\bel
Let $\cG=\cL,\cR,\cL\cR,\cC,\cK,$ then $exp,ln: T_{\cG^{(j)}} \rightleftarrows \cG^{(j)}$ for $j\ge1.$
\eel
The exponential/logarithmic maps of filtration-unipotent operators are well known, e.g.
  \S III.5.4 of \cite{Bourbaki.Lie} or pg.27 of \cite{Serre.Lie}.
The lemma is well-known for $\MapX,$ with the filtration $I^\bullet\cdot \RmX,$ and the groups $\cR,\cK,$ see e.g. \cite[\S3.5]{Belitskii-Kerner} and \cite[Theorem 3.6]{BGK.20}. We need this result for singular targets, $\Maps.$ Therefore we just add the needed ingredients in our more general situation.
 \bpr
 It is enough to verify the statement for the groups $\cL^{(j)},$ $\cR^{(j)},$ $\cC^{(j)}.$ Then the case of composite groups, $(\cL\cR)^{(j)},$ $\cK^{(j)}$ follows.

$\pmb{exp}\!: T_{\cG^{(j)}}\to\cG^{(j)}.$ For every $\xi \in T_{\cG^{(j)}}$ the action of the  operator $exp(\xi)\circlearrowright \RmX$ is well defined. Thus it is enough to verify: $exp(\xi)$  acts on the subset $\Maps $ and is unipotent for the filtration $\tM^\bullet\sset\RmX.$

 The part $exp(\xi)\circlearrowright\Maps$ is verified case-wise. Let $f\in \Maps,$ then:
 \bei
 \item $\cR\!:$ \quad $exp(\xi_X)\in Aut_X,$ therefore  $exp(\xi_X)(f(x))=f(exp(\xi_X)(x))\in \Maps.$
 \item $\cL\!:$ \quad $exp(\xi_Y)\in Aut_Y,$ therefore   $exp(\xi_Y)(f)\in \Maps.$
 \item $\cC\!:$ \quad here $\xi(J_Y)\sseteq R_{X\times Y}\cdot J_Y=0\sset R_{X\times Y}.$ Denoting
 $\tf\!:=exp(\xi)f $ we get: $\tf^\sharp(J_Y)\sseteq R_{X\times Y}\cdot J_Y,$ i.e. $\tf^\sharp\sin Hom(R_Y,R_X),$ i.e. $\tf\in \Maps.$
 \eei
Finally, the $\tM^\bullet$-unipotence of $exp(\xi)$ is immediate, $(exp(\xi)-Id)\tM^d=(\sum_{l\ge1}\frac{\xi^l}{l!})\tM^d\sseteq \tM^{d+j}.$

\medskip

$\pmb{ln}\!: \cG^{(j)}\to T_{\cG^{(j)}}.$  For every $g \in  \cG^{(j)} $ the action of the  operator $ln(g)\circlearrowright \RmX$ is well defined. Thus it is enough to verify: $ln(g)\in T_\cG$ and is nilpotent for the filtration $\tM^\bullet\sset\RmX.$

  \bei
 \item   The facts $ln(g)\in Der_X$ for $g\in Aut_X^{(1)},$ and $ln(g)\in Der_Y$  for $g\in Aut_Y^{(1)}$ are well known, see e.g.
   \cite[\S3.5]{Belitskii-Kerner} and \cite[Theorem 3.6]{BGK.20}.

 \item  For the group $\cC<Aut_{X\times Y}$ we denote $\xi:=ln(g).$ Thus $\xi\in Der_{X\times Y},$ as before. Moreover,
  we have $\xi(R_{X\times Y})\sseteq R_{X\times Y}\cdot \cm_Y,$ as $g\cdot \cm_Y\sset R_{X\times Y}\cdot \cm_Y.$ And $\xi$ is $x$-linear, as $g$ preserves the fibres of the projection $X\times Y\to Y.$ Therefore $\xi\in R_{X\times Y}\otimes_{R_Y} Der_Y Log(\cm_Y).$
 \eei
 Finally, the $\tM^\bullet$-nilpotence of $ln(g)$ is immediate, $ln(g)\tM^d=-(\sum_{l\ge1} \frac{(1-g)^l}{l})\tM^d\sseteq \tM^{d+j}.$
 \epr

\section{The relevant Artin approximation  and the $T_\cG,\cG$-versions of the Artin-Rees lemma}

\subsection{Maps of finite singularity type}\label{Sec.K.finite.maps}
This notion is used below only to quote the Artin approximation   for the particular case of non-finitely determined maps, over the analytic rings $\quots{\k\{x\}}{J_X},\quots{\k\{y\}}{J_Y},$ for the groups $\cG=\cL,\cL\cR.$

 A map $X\stackrel{f}{\to}Y$ is called {\em of finite singularity type} if its restriction to the critical locus,
  $X\supset Crit\stackrel{f_|}{\to}Y,$ is a finite morphism.
   In the classical case of $\MapC$ this notion is classical,  \cite[pg.224]{Mond-Nuno}.
For the basic properties  in the case $\MapX$ see
   \cite[\S3.6.2]{Kerner.Group.Orbits}, for the case $\Maps$  see \cite[\S4.4]{Kerner.LRAP}.

\medskip

\noindent When $char(\k)>0$ and [$J_X\!\neq\!0$ or $J_Y\!\neq\!0$] we assume in addition the following ``lifting properties":
\bei
\item   the approximate integrability of derivations, \cite[Definition 4.1]{Kerner.LRAP},
\item and the approximate lifting of automorphisms
   \cite[equation (16)]{Kerner.LRAP}.
 \eei
We do not restate them in detail, as they are not used below.

\subsection{Artin approximation for $\cG$-action}\label{Sec.Background.AP}
Let $\cG$ be one of the groups $\cL,\cR,\cL\cR,\cC,\cK.$ Fix two map-germs $X\stackrel{f,\tf}{\to}Y.$ We want to establish their $\cG$-equivalence, $\tf\stackrel{\cG}{\sim}f,$ i.e. to prove that $\tf$ belongs to the orbit $\cG f.$

\subsubsection{Ordinary $\cG$-Artin approximation (``from formal equivalence to Nash/analytic equivalence")}\label{Sec.Background.AP.GAP} The verification of the $\cG$-equivalence often goes at the level of formal power series. Namely, we work with the completions of the rings at their maximal ideals: $\hR_X=\quots{\hk[\![x]\!]}{J_X},$ $\hR_Y=\quots{\hk[\![y]\!]}{J_Y}.$
 (If $\k$ is   a field then  $\hk=\k.$)
 Then $\hat\cG$ is one of the corresponding groups:
\beq
Aut_\hk(\hR_X),\hspace{1.2cm} Aut_\hk(\hR_Y),\hspace{1.2cm}  Aut_\hk(\hR_Y)\times Aut_\hk(\hR_X), \hspace{1.2cm} \hat\cK:= \hat\cC\rtimes  Aut_\hk(\hR_X).
\eeq
Having established the formal equivalence, ``$\tf\stackrel{\hat\cG}{\sim}f$", we want to deduce the ordinary (i.e. Nash/analytic) equivalence  $f\stackrel{\cG}{\sim}\tf.$  This goes via the Artin approximation.
 In the case [$\cG=\cL,\cL\cR$ and $f$ is not finitely-$\cG$-determined] we assume:
 \beq\label{Eq.A.P.Assumptions.for.L,LR,groups}
\text{either $R_X=\quots{\k\bl x\br}{J_X} $   or [$R_X=\quots{\k\{ x\}}{J_X}$ and $f$ is of finite singularity type].}
\eeq
\bthe
The formal equivalence, $\tf\stackrel{\hat\cG^{(j)}}{\sim}f,$ is approximated by the Nash, resp. analytic  equivalence, $\tf\stackrel{\cG^{(j)}}{\sim}f.$
\ethe
Namely, if $\tf=\hg\cdot f $ for some $\hg\sin \hat\cG^{(j)},$ then for each $d\sin \N$ there exists $g_d\sin \cG^{(j)}$ satisfying:
 $\tf=g_d f$ and $\hg \equiv g_d\ mod\ \cm^d.$

\bei
\item
For the groups $\cR,\cK$ and $\k$-smooth target, $Y=(\k^m,o),$ this statement is an immediate application of the classical Artin approximation, see e.g. \cite[\S2.3.v]{Kerner.Group.Orbits}.
\item For the group $\cL\cR,$  and for $\cK$ when $Y\ssetneq(\k^m,o),$  this is a non-trivial statement\footnote{In fact, the $\cL,\cL\cR$ properties fail for analytic maps that are not of finite singularity type, due to Osgood-Gabrielov-Shiota counterexample, \cite{Shiota.1998}.}. The cases $\cL,\cL\cR$ are proved in \cite{Kerner.LRAP}, Theorem 3.4 for Nash maps, and Theorem 5.1 for analytic maps that are of finite singularity type. The general case $\cG^{(j)}$ (in particular for $\cG=\cL\cR,\cK$) is done in Theorem 4.3 (and \S4.4) of \cite{Kerner.G-AP}.
\eei

\subsubsection{Strong $\cG$-Artin approximation (``from order-by-order equivalence to ordinary equivalence")}\label{Sec.Background.AP.SAP} In some cases we cannot even reach the conclusion on formal power series  ``$\tf\stackrel{\hat\cG}{\sim}f$".
 We only get the equality modulo higher order terms:
\beq
  ``\tf \equiv g_d\cdot f\ mod\ \cm^d \quad  \text{ for each $d>1$ and a corresponding } g_d\sin \cG."
   \eeq
   Here the sequence of group-elements $g_\bullet$ is obtained very indirectly. Often we cannot force the sequence $g_\bullet$  to converge in the $\cm^\bullet$-adic topology. Thus we cannot just conclude
  ``$\tf=(\liml_{d\to\infty} g_d) f$".

  Restate this via the filtration-topology. Take the orbit $\cG f$ and its $\cm$-adic neighborhoods $\cG f+\cm^\bullet\scdot \RmX \sseteq \RmX.$
   Suppose $\tf\sin \cG f+\cm^d\scdot \RmX $ for each $d\sin \N.$ Then $\tf\sin \cap_d (\cG f+\cm^d\scdot \RmX ),$ the latter being the filtration closure
    of the orbit $\cG f\sseteq  \RmX.$ The  Strong Artin approximation ensures that the orbit is closed, in the following way.

  Let $R_X,R_Y$ be as in \S\ref{Sec.Background.Germs.Rings}, and take a map $X\stackrel{f}{\to}Y.$
   Take a filtration $\tM^\bullet$ equivalent to $\cm^\bullet\cdot \RmX.$ Take the corresponding subgroups $\cG^{(\bullet)}<\cG.$
     In the case [$\cG=\cL,\cL\cR$ and the map is formal or analytic, and not finitely $\cG$-determined]     we assume:
\beq\label{Eq.SAP.Assumptions.for.L.LR.group}
\text{the base field is $\aleph_0$-complete, and, in the analytic case,   $f$ is of finite singularity type
.}
  \eeq
Examples of   $\aleph_0$-complete fields are: finite fields and uncountable  fields that are algebraically closed. The fields $\Q,\R$ are not
    $\aleph_0$-complete,  see Theorem 5 in \cite{Popescu.Rond}.

  \bthe\cite[Theorem 4.3 and Proposition 4.5]{Kerner.G-AP}\label{Thm.SAP.G} The orbits $\cG f$ and $\cG^{(j)}f$  are closed for the filtration $\cm^\bullet\cdot \RmX.$
  \ethe
  Namely: $\cap_d (\cG  f+\cm^d\cdot \RmX)=\cG f $ and
   $\cap_d (\cG^{(j)} f+\cm^d\cdot \RmX)=\cG^{(j)} f,$  for each $j\ge0$ and   each $f\in \Maps.$
 Equivalently, if $\tf\not\in\cG^{(j)} f,$ then $\tf\not\in\cG^{(j)} f+ \cm^d\cdot \RmX$ for some $d\sin \N.$

   Yet more explicitly,  for any maps $X\stackrel{\tf,f}{\to}Y$ there exists the Artin function $\be_j:=\be(j,f,\tf,\cG,\tM^\bullet)$ that ensures:
    \beq\label{Eq.Artin.Function}
    \text{If } \quad \tf\sin \cG^{(j)} f+ \cm^{\be_j}\cdot \RmX  \quad  \text{  then } \quad \tf\sin \cG^{(j)} f.
    \eeq
\beR\label{Rem.Artin.function.invariant.on.G.orbits}
This Artin function $\be_j$ depends only on the $\cG^{(0)}$-orbit of the pair $f,\tf,$ i.e.
 \beq
 \be(j,f,\tf,\cG,\tM^\bullet)=\be(j,g\scdot f,g\scdot \tf,\cG,\tM^\bullet) \quad \text{ for each }\quad  g\in \cG^{(0)}.
 \eeq
  Indeed, suppose
  $g\tf\in \cG^{(j)}gf+\cm^{\be_j}\scdot \RmX$ for some $g\in \cG^{(0)}.$ Then  $ \tf\in g^{-1}\cG^{(j)}gf+\cm^{\be_j}\scdot \RmX.$
   By the normality $\cG^{(j)}\lhd \cG^{(0)},$ see Remark \ref{Rem.on.FIltrations.of.G}.iii, we get
    $ \tf\in  \cG^{(j)} f+\cm^{\be_j}\scdot \RmX.$ Hence $ \tf\in  \cG^{(j)} f .$ Therefore  $ g\tf\in  g\cG^{(j)} f = \cG^{(j)} gf.$

Recall that for the filtration $\tM^\bullet=\cm^\bullet\cdot\RmX$  we have $\cG^{(o)}=\cG.$
\eeR

\subsection{The $T_\cG f$-version of the Artin-Rees lemma}\label{Sec.Artin-Rees.TG}
Fix a filtration $\tM^\bullet\sset \RmX$ equivalent to $I^\bullet\cdot \RmX,$ see \S\ref{Sec.Background.Filtrations.of.Maps}.
 In the proof of Theorem \ref{Thm.Equivalence.Change.base.ring} we use the following statement:
\beq\qquad
T_\cG f\cap \tM^{d_j}\sseteq   T_{\cG^{(j)}}f \qquad
\text{ for each }j\!\in\! \N \text{ and a corresponding }1\ll d_j\!\in\! \N.
\eeq
For $\cG=\cR,\cC,\cK$ this follows straight from the   Artin-Rees lemma, \cite[\S8]{Matsumura}. Indeed, in this case  the image tangent space, $T_\cG f\sseteq \RmX,$ is a finitely-generated $R_X$-submodule, \S\ref{Sec.TG.action.on.Maps}.

For $\cG=\cL,\cL\cR$ the image tangent space  $T_\cG f\sseteq \RmX $ is only an $R_Y$-submodule,  not an $R_X$-submodule.
 And $T_{\cL\cR}f$ is not necessarily finitely-generated over $R_Y.$   Yet we have:
\bel
Let $f\sin \Maps  \cap I\cdot \RmX.$ Then $T_{\cL }f\cap \tM^{d_j}\sseteq T_{ \cL ^{(j)}}f $ and $T_{\cL\cR}f\cap \tM^{d_j}\sseteq T_{(\cL\cR)^{(j)}}f,$ for each $j\sin \N$ and a corresponding $1\ll d_j\sin \N.$
\eel
Here $f$ is not necessarily finitely-$\cG$-determined.

\noindent This statement seems to be known  for $\cG=\cL\cR$ in the particular case of $\MapC. $   For $\MapX$ it is proved in \cite[Lemma 6.11]{Kerner.Group.Orbits}.
Essentially the same proof holds in our (more general) case, we repeat it for completeness for $\cG=\cL\cR.$
\bpr
By our assumption the filtrations  $\tM^\bullet$ and $I^\bullet\cdot \RmX$ are equivalent. Therefore we can assume  $\tM^d=I^d\cdot \RmX.$
 Consider the quotient of the image tangent spaces $\quots{T_{\cL\cR}f}{T_\cR f}.$ This is a finitely-generated $R_Y$-module, e.g. by the surjection $T_\cL f\twoheadrightarrow \quots{T_{\cL\cR}f}{T_\cR f}.$
  It acquires the filtration $\quots{T_{\cL\cR}f\cap I^\bullet\cdot \RmX+T_\cR f}{T_\cR f}\sset \quots{T_{\cL\cR}f}{T_\cR f},$ by $R_Y$-submodules.
 These submodules are finitely generated as well, by  Remark \ref{Rem.Noetherian.Submodules}.

  This filtration is separated, i.e. $\cap_\bullet\big(\quots{T_{\cL\cR}f\cap I^\bullet\cdot \RmX+T_\cR f}{T_\cR f}\big)=0.$
   Therefore
   \beq
   \quot{T_{\cL\cR}f\!\cap\! I^{d_j}\cdot \RmX\!+\!T_\cR f}{T_\cR f}\sset  (f^\sharp(y)+\cm_\k\cap I)^{j+1} \!\cdot\! \quot{T_{\cL\cR}f}{T_\cR f},
     \eeq
    for each $j\!\in\! \N $   and a corresponding $j\!\ll\! d_j\!\in\! \N.$
\\
Thus
     $T_{\cL\cR}f\cap I^{d_j}\cdot \RmX \sset (f^\sharp(y)+\cm_\k\cap I)^{j+1}\cdot T_{\cL\cR}f+T_\cR f.$
      As $f\sin I\cdot\RmX,$ we  get
 $ T_{\cL\cR}f\cap I^{d_j}\cdot \RmX \sset    T_{(\cL\cR)^{(j)}}f+ T_\cR f.$
 Finally, comparing the two sides in this last inclusion  we get
\beq
T_{\cL\cR}f\cap I^{d_j}\cdot \RmX \sset  T_{(\cL\cR)^{(j)}}f+ T_\cR f\cap I^{j }\RmX.
\eeq
 And then  $ T_\cR f\!\cap\! I^{j }\RmX\sseteq  T_{\cR^{(j-l)}} f,$  for $j\!\gg\!1$ and a constant $l\sin \N,$ by the classical Artin-Rees.
 \epr

\subsection{The $\cG f$-version of the Artin-Rees lemma}\label{Sec.Artin-Rees.G}   
 Fix a ($\k$-Nash, resp. $\k$-analytic, formal) map $f\sin \Maps,$ and let $\cG$ be one of the groups $\cL,\cR,\cL\cR,\cC,\cK.$
  Here $f$ is not necessarily finitely-$\cG$-determined.
 In the case [$\cG=\cL,\cL\cR$ and $f$ is formal or analytic]  we take the assumptions \eqref{Eq.SAP.Assumptions.for.L.LR.group}.
  Take the filtration $\tM^\bullet=\cm^\bullet_X\scdot\RmX.$
 Take the Artin function $\be$ of  \eqref{Eq.Artin.Function}  for the equation  [$f=g f,$ $g\in \cG$].
\bel
$\cG f\cap (\{f\}+\cm^{\be_j}\cdot \RmX)\sset \cG^{(j)}f$ for each $j\ge1.$
\eel
\bpr
Take any $\tf\in \cG f\cap (\{f\}+\cm^{\be_j}\cdot \RmX).$ Thus $\tf=\tg_o\cdot f $ for some $\tg_o\in \cG.$ Consider this $\tg_o$ as an approximate solution to the equation $[gf=f,\ g\in \cG].$ By Theorem \ref{Thm.SAP.G} we get an ordinary solution, $g_o f=f,$ $g_o\in \cG,$
 which moreover satisfies $g_o\equiv \tg_o\ mod\ \cm^j .$
  Thus $\tg_o g^{-1}_o\in \cG^{(j)}.$ Altogether:
 \beq
\hspace{3cm} \tf=\tg_o f=(\tg_o g^{-1}_o)g_o f=(\tg_o g^{-1}_o)  f\in \cG^{(j)} f.
 \hspace{4cm}\epr
 \eeq

\bex\label{Ex.G.-fin.det.iff.Gjj.fini.det}
As an immediate application, the following conditions are equivalent:
\bei
\item $f$ is $\cG$-finitely determined
\item  $f$ is $\cG^{(j)}$-finitely determined for some $j\sin \N $
 \item  $f$ is $\cG^{(j)}$-finitely determined for every $j\sin \N.$
\eei
\eex

\beR
 For particular applications one needs a quantitative version of this lemma. No reasonable bounds on the Artin function $\be_j$ seem to be known, as the condition $gf=f$ is a highly non-linear system of equations. However, if $f$ is finitely $\cG^{(j)}$-determined, then a trivial bound comes from
  $\cG^{(j)}f\supseteq (\{f\}+\tM^d)\cap \Maps.$ Here $d$=the order of determinacy.
\eeR

\section{$\cG_\k$-equivalence of maps vs their $\cG_\K$-equivalence}\label{Sec.K.vs.k}
Let $R_X$ be one of the rings $\quots{\k[\![x]\!]}{J_X},   \quots{\k\{x\}}{J_X},  \quots{\k\bl x\br}{J_X},$ as in \S\ref{Sec.Background.Germs.Rings}.
 Take any  ring extension $\k\hookrightarrow \K $.
 Accordingly take the   extended rings   $R_{X, \K }:=\quots{ \K [\![x]\!]}{(J_X)},$   $\quots{ \K \{x\}}{(J_X)},  \quots{ \K \bl x\br}{(J_X)}$.
 We get the groups
 \beq
  \cR_ {\K }=Aut_ \K(R_{X, \K }),\quad
 \cL_ {\K }=Aut_ \K(R_{Y, \K }),\quad
 \cC_\K< Aut_\K(R_{X\times Y,\K}),\quad
 \cK_\K=\cC_\K\rtimes \cR_\K.
 \eeq
Every element $g\in Aut_{X\times Y}$ is $\k$-linear, thus it extends (uniquely) to a $\K$-linear element $g_\K \in Aut_{X_\K\times Y_\K}.$
 This gives the group-extension $\cG \le \cG_ {\K }.$

 Every map of $\k$-germs $X\stackrel{f}{\to}Y$ extends (uniquely) to a map  $X_\K\stackrel{f_\K}{\to}Y_\K.$
 We compare the $\cG$-equivalence of $\k$-maps to their $\cG_ {\K }$-equivalence.

\subsection{}\label{Sec.K.vs.k.main.Theorem} As in \S\ref{Sec.G.filtered} we fix a filtration $\tM^\bullet\sset\RmX $ equivalent to the filtration $I^\bullet\cdot \RmX.$
 Accordingly we define $\tM_\K:=R_{X,\K}\cdot \tM^\bullet\sset R^{\oplus m}_{X,\K}.$
 This gives the filtration on the space of maps,
 \beq
 M^\bullet:=\Maps\cap \tM^\bullet\quad \text{ and similarly }\quad M^\bullet_\K:=\MapsK\cap  \tM^\bullet_\K.
 \eeq
As in \S\ref{Sec.G.filtered}  we get the subgroups $\cG^{(j)}\le \cG$ and  $\cG^{(j)}_\K\le \cG_\K.$
 Below we prove:
 \beq
 \cG^{(j)}_{\K }f\cap \Maps=\cG^{(j)} f.
 \eeq

   We assume: $f\sin I\cdot \RmX.$
 It is not necessarily finitely-$\cG$-determined.
    In the case [$\cG=\cL,\cL  \cR$ and $R_X=\quots{\k\{ x\}}{J_X}$ and $f$   is not finitely-$\cG$-determined]   we also assume that $f$ is of finite singularity type, see \S\ref{Sec.K.finite.maps}.
\bthe\label{Thm.Equivalence.Change.base.ring}
Take two (formal, resp. $\k$-analytic, $\k$-Nash) maps $X\stackrel{f,\tf}{\to}Y$.
       \bee[\bf 1.]
\item Suppose   $char(\k)\= 0$ and the ring extension $\k\sset \K$ is faithfully-flat, \S\ref{Sec.Background.FF.extensions}.
Let $j\!\ge\!1$.
$$
\text{ If  }f\!\stackrel{\cG^{(j)}_{\K }}{\sim}\!\tf \ \text{ then }\ f\!\stackrel{\cG^{(j)}}{\sim}\!\tf.
$$
\item Suppose the base field ($\k,$ resp. $\k_o$)   is  algebraically   closed. In the case [$\cG=\cL,\cL  \cR$ and $R_X$ is analytic or formal, and $f$ is not finitely-$\cG$-determined]  we also assume that the base field is $\aleph_0$-complete, see \S\ref{Sec.Background.AP.SAP}.

   If $f\stackrel{\cG_{\K }}{\sim}\tf$ then $f\stackrel{\cG}{\sim}\tf$.   Moreover, for $I=\cm_X $ and  any $j\ge0:$ \quad
  if $f\stackrel{\cG^{(j)}_{\K }}{\sim}\tf$ then $f\stackrel{\cG^{(j)}}{\sim}\tf$.
\eee\ethe
\bpr It is enough to prove both statements in the formal case, $R_X=\quots{\hk[\![x]\!]}{J_X}$. In the analytic case, $\quots{\k\{x\}}{J_X},$ and in the  Nash  case, $\quots{\k\bl x\br}{J_X},$   one invokes the $\cG$-Artin approximation,
 \S\ref{Sec.Background.AP.GAP}.

 \medskip

The proof of part 1  is based on  the unipotence. We reduce the question to the level of tangent spaces. And there we have ``a linear problem", manageable by Commutative Algebra.
 (With exactness of the functor $\otimes \K$ and the $T_\cG$-Artin-Rees lemma.)

For part 2 we translate the conditions into a countable system of polynomial equations (over the base field $\k=\bk,$ resp. $\k_o=\bk_o$). Each finite subsystem is resolvable. Thus we get: $\tf\in \cG^{(j)} f+\cm^d\cdot \RmX$ for all $d\in \N.$
 Finally the Strong Artin approximation is invoked.

\bee[\bf 1.]
\item  Suppose $\tf=g_\K f$ for an element $g_{\K }\sin \cG^{(j)}_{\K }$, with $j\ge1$. Present $g_{\K }$ as the exponential of a tangent vector,
 $g_{\K }=exp(\xi^{\K }_d),$ see   \S\ref{Sec.TG.to.G.exp}.  Here  $\xi^{\K }_d\sin T_{\cG^{(j)}_{\K }},$   while the index $d$ is defined by
 the $\tM^\bullet_\K$-order, $d:=ord (\xi^\K_d f)-ord(f)\ge j,$ see \S\ref{Sec.Background.Filtrations.of.Maps}.
\\
Note that $T_{\cG^{(j)}_{\K }}\= \K \otimes T_{\cG^{(j)}}.$ Here (and below) the tensor product $\otimes $   is taken over $\k,$ and  if $\K$ is
 not a finitely-generated $\k$-module, then one takes the completed version $\hotimes$ \wrt   the filtration $\cm^\bullet_X\scdot \RmX$.
\bei
\item
  We claim: $\xi^{\K }_d f\sin T_{\cG^{(j)}}f+ \tM^{ord(f)+d+j}_\K .$ 
 Indeed,  $ord (\xi^{\K }_d  f) = ord(f) +d.$
  Therefore:
\beq
 exp(\xi^{\K }_d)f-f-\xi^{\K }_d f= (\xi^{\K }_d)^2(\dots)\cdot f\sin \tM^{ord(f)+d+j}_\K  .
\eeq
   As $g_\K f,f\sin \RmX$, we get
   $\xi^{\K }_d f\sin \RmX+ \tM^{ord(f)+d+j}_\K.$ 

    Pass to the quotient by the $\K$-submodule $\tM^{ord(f)+d+j}_\K\sset \RmXK$ to get:

   \beq
 [\RmX] \ni \  [\xi^{\K }_d f] \ \sin    [\K\otimes T_{\cG^{(j)}}f]=\K\otimes[ T_{\cG^{(j)}}f].
\eeq

  Thus $[\xi^{\K }_d f]\sin [ \K \otimes T_{\cG^{(j)}}f]\cap [\RmX].$

   By faithful-flatness  we have  $[ \K \otimes T_{\cG^{(j)}}f]\cap [\RmX]=[  T_{\cG^{(j)}}f] ,$
    Lemma  \ref{Thm.Faithful.Flatness}.
 Thus
 $\xi^{\K }_d f\sin T_{\cG^{(j)}}f+ \tM^{ord(f)+d+j}_\K $.

\medskip

\item
 Therefore we can expand
 $\xi^{\K }_d=\xi_d+\txi^{\K }_{d+j}$ and present $\tf=exp(\xi_d+\txi^{\K }_{d+j})f.$
 Here
 \beq
 \xi_d\sin T_{\cG^{(j)}}, \quad   \txi^{\K }_{d+j}\sin  \K \otimes T_{\cG^{(j)}}
 \quad
 ord(\xi_d f)=   ord( f)+d, \quad    ord(\txi^{\K }_{d+j} f)\ge   ord(f)+d+j.
  \eeq
  Now we have:  $exp(-\xi_d)\tf=exp(\xi^{\K }_{d+j})f,$ using the Baker-Campbel-Hausdorff formula,   \cite[\S IV.7]{Serre.Lie}.
 Here $\xi^{\K }_{d+j}\sin  \K \otimes T_{\cG^{(j)}}.$  By the direct check: $ord \xi^{\K }_{d+j} (f)\ge ord(f) +d+j.$

 Replace $\tf$ by $exp(-\xi_d)\tf$ and iterate the whole argument for $d\rightsquigarrow d+j$.
 Thus, starting from $\xi^{\K }_d,$ we get elements $\xi_d,\xi_{d+j},\xi_{d+2j},\cdots\sin  T_{\cG^{(j)}}$ satisfying:
 \beq
 exp(-\xi_{d+lj})\cdots  exp(-\xi_{d+j})\cdot  exp(-\xi_{d})\tf=exp(\xi^{\K }_{d+(l+1)j})f.
 \eeq
 Here
 $ord (\xi^{\K }_{d+(l+1)j}f)\ge ord(f)+d+(l+1)j.$

\item
We want to take the  product of the constructed elements, $g:=\lim_{l\to\infty} (e^{-\xi _{d+lj}}\cdots e^{-\xi _d}).$
 This limit does not necessarily exist, as the vectors $\xi_\bullet$ do not converge in whichever sense.
  However, observe that   $\xi_d f\sin T_{\cG^{(j)}}f\cap M^d $ for each $d.$  Apply the Artin-Rees  property, \S\ref{Sec.Artin-Rees.TG}.
   We get:  $\xi_d f\sin T_{\cG^{(i_d)}}f$, with $\lim_d i_d=\infty$.
  Thus we rerun the proof, getting the adjusted vectors $\ze_d\sin T_{\cG^{(i_d)}}.$
   And for these vectors the (formal) limit $\hg:=\lim_{l\to\infty} (e^{-\ze _{d+lj}}\cdots e^{-\ze _d})\sin \cG^{(j)}$ does exist.

\eei

{\hspace{-1cm}Altogether we have $\tf=\hg\cdot f.$  This proves the first statement in the formal case, $R_X=\quots{\k[\![x]\!]}{J_X}$.}

\

\item (We prove only the case of $\cL\cR,$ the other cases are similar and simpler.)
 Recall, $R_X=\quots{\hk[\![x]\!]}{J_X},$ $R_Y=\quots{\hk[\![x]\!]}{J_Y}$.
We should resolve the condition:
\beq\label{Eq.inside.proof.K.vs.k}
\Phi_Y\circ f=\tf\circ\Phi_X, \hspace{1cm} \text{for} \hspace{1cm} \Phi_Y\sin \cL   \hspace{0.5cm} \text{and} \hspace{0.5cm} \Phi_X\sin \cR .
\eeq
\bei
\item We have $\hk=\quots{\k_o[\![t]\!]}{\ca},$ for an algebraically closed field $\k_o.$
 Take the Taylor expansions,
 \beq\label{Eq.Expanding.Autos}
 \Phi_X(x)=\sum_{|m|\ge1,|l|\ge0} C^{(x)}_{m,l}\cdot t^l\cdot x^m \quad\text{ and  }\quad \Phi_Y(y)=\sum_{|m|\ge1,|l|\ge0} C^{(y)}_{m,l}\cdot t^l\cdot y^m.
 \eeq
  Here $m,l$ are multi-indices,
  $x^m,y^m,t^l$ are monomials, and
  $\{C^{(x)}_{m,l}\}$,   $\{C^{(y)}_{m,l}\}$ are unknowns in $\k_o$. By comparing the coefficients of the $x,t$-monomials, the condition
    $\Phi_Y\circ f=\tf\circ\Phi_X,$ together with $\Phi_X(J_X)=J_X\sset \k[\![x]\!]$ and $\Phi_Y(J_Y)=J_Y\sset \k[\![y]\!],$  is transformed into a countable system of polynomials equations, $\{P_l(\{C^{(x)}_{**}\},\{C^{(y)}_{**}\})=0\}_l$.
     Each polynomial here is in a finite number of variables.

     To these equations we add the conditions: $\Phi_X,\Phi_Y$ are invertible. Which means: the $\k_o$-matrices $\{C^{(x)}_{m,0}\}_{|m|=1},$
       $\{C^{(y)}_{m,0}\}_{|m|=1}$ (of finite size) are non-degenerate. These non-degeneracy conditions are presented as polynomial equations
       $det[\{C^{(x)}_{m,0}\}_{|m|=1}]\cdot z_x=1,$ $det[\{C^{(y)}_{m,0}\}_{|m|=1} ]\cdot z_y=1,$ with new unknowns $z_x,z_y.$

      Take the ideal generated by all these polynomials, $\cP\sset \k_o[C^{(x)}_{**},C^{(y)}_{**},z_x,z_y].$
 By our assumption     the corresponding system of polynomial equations is solvable over the ring $ \K $.

We claim: {\em each finite subsystem of this system  is solvable over $\k_o$.}
 Indeed, this finite subsystem       defines an algebraic subscheme in a  (finite dimensional) affine space over (an algebraically closed field) $\k_o.$
       If this finite subsystem is not solvable, then this subscheme has no points over $\k_o$.
       By Hilbert Nullstellensatz we get: the defining ideal of this subscheme is the whole polynomial ring of this affine space.
        But then the ideal $\cP$ contains $1$, hence         the whole system cannot have
         solutions over any ring $ \K $.

\item (The case of the whole group $\cL\cR.$) Fix an integer $N\gg1.$ From the infinite set of elements of $\cP$ we extract  those containing only unknowns of bounded order:
\beq
\cP\cap  \k_o[\{C^{(x)}_{m,l}\}_{|m|+|l|\le N},\{C^{(y)}_{m,l}\}_{|m|+|l|\le N},z_x,z_y].
\eeq
      This is an ideal in the Noetherian ring, hence is finitely generated. Thus the corresponding system of (polynomial) equations is  solvable over $\k_o.$

    Let $\Phi_{X,N}(x),$ $\Phi_{Y,N}(y)$ be  power series   corresponding to this solution.  (We assign arbitrary values to the constants   $\{C^{(x)}_{m,l}\}_{|m|+|l|> N},\{C^{(y)}_{m,l}\}_{|m|+|l|> N},$ e.g. zeros.)
     They define the automorphisms $\Phi_{X,N}\circlearrowright(\k^n,o),$ $\Phi_{Y,N}\circlearrowright(\k^m,o).$ But these automorphisms do not necessarily restrict to automorphisms of $X,Y,$ as $\Phi_{X,N},\Phi_{Y,N}$ preserve the ideals $J_X,J_Y$ only  up to terms of order $N.$
      Yet (by the strong Artin approximation for $Aut_X,Aut_Y$), we can adjust $\Phi_{X,N},\Phi_{Y,N}$ by  terms of order $d_N\gg1,$ so that the adjusted versions do define automorphisms of $X,Y,$ thus give elements of $\cR,\cL.$ Altogether, we get an element $g \sin \cL  \cR$ satisfying: $g  \tf\in f+\cm^d\cdot \RmX.$ Therefore $\tf\in \cG f+ \cm^d\cdot \RmX,$ where $d\gg1,$ as $N\gg1.$
 Moreover, for any $d\in \N$ we can rerun the proof (with $N_d\gg1$) to ensure $\tf\in \cG f+ \cm^d\cdot \RmX.$

Now invoke the strong-$\cG$-Artin approximation, Theorem \ref{Thm.SAP.G}, to get $\tf\stackrel{\cG}{\sim}f.$

\item  The case of the groups $(\cL\cR)^{(j)}$ for $j\ge0$ is done by a slight modification of the argument. We are looking for an element
 $(\Phi_Y,\Phi_X)\sin (\cL\cR)^{(j)}$ satisfying conditions \eqref{Eq.inside.proof.K.vs.k}.
   We expand $\Phi_Y,\Phi_X$ as in \eqref{Eq.Expanding.Autos}. Thus  conditions \eqref{Eq.inside.proof.K.vs.k} are
       converted  into  equations in $C^{(x)}_{**},C^{(y)}_{**}.$ These equations are polynomial. As before, we get an approximate solution
       $g=(\Phi_Y,\Phi_X),$ i.e.
\beq
\Phi_Y\tf\stackrel{mod\ \cm^d}{\equiv}\Phi_X f, \quad \Phi_Y\stackrel{mod\ \cm^d}{\in}  \cL ^{(j)}, \quad
 \Phi_X\stackrel{mod\ \cm^d}{\in}  \cR ^{(j)}.
\eeq
By approximating $g$  we can assume: $g_d\sin (\cL\cR)^{(j)}.$
  Finally invoke the Strong Artin approximation again,   to get $g\sin (\cL\cR)^{(j)}$ satisfying $g\tf=f.$
     \epr
\eei
\eee

\bex
We restate part 2 in the most important case. Take any algebraically closed field $\k.$
 Take two (Nash/analytic/formal) map germs $(\k^n,o)\stackrel{f,\tf}{\to}(\k^m,o).$
    Let $\cG$ be one of the groups $\cL,\cR,\cL\cR,\cC,$ $\cK.$
    If [$R_X$ is analytic or formal, and $\cG\!=\!\cL,\cL\cR,$ and   $f$ is not finitely-determined] then we assume in addition
     [$\k$ is  uncountable and   $f$ is of finite singularity type].
\[
\text{  If }f\stackrel{\cG_\K}{\sim}\tf \text{ then  } f\stackrel{\cG}{\sim}\tf.
\]
\eex
\beR
The proof of part 2 of Theorem \ref{Thm.Equivalence.Change.base.ring} gives more than just $\tf\stackrel{\cG}{\sim}f.$ In this proof the $\cG_\K$-equivalence is ``approximated" by $\cG$-equivalence in the following sense. Suppose
  $\tf=g_\K\cdot f,$ where $g_\K\in \cG_\K$ is presented by power series in $x,y,$ as in \eqref{Eq.Expanding.Autos}. Suppose the coefficients $\{C^{(x,\K)}_{**}\},$  $\{C^{(y,\K)}_{**}\}$ satisfy a countable collection of polynomial equations over $\k_o.$
     Then
  $\tf=g\cdot f,$ where $g\in \cG$ and the coefficients  $\{C^{(x)}_{**}\},$  $\{C^{(y)}_{**}\}$ of $g$ satisfy the same equations.
 Thus if $g_\K$ satisfies some functional equations (not necessarily of implicit function type) or differential equations, all over $\k_o,$ then
  $g$ will satisfy the same equations.
\eeR

\subsubsection{Remarks on the assumptions of the first statement in Theorem \ref{Thm.Equivalence.Change.base.ring}}\label{Sec.Remarks.Sharpness} \

 They are   annoying, but (it seems) cannot be weakened.
\bee[\bf i.]
 \item\ [The case $j=0$] In any characteristic, suppose a field $\k$ is not algebraically closed, and  $c\sin \k\smin \k^p,$ i.e. $c$ is not a $p$'th root.
 Then $c\cdot x^p\stackrel{\cR^{(0)}}{\sim}x^p$ over $\bk,$ but not over $\k.$

  Similarly, $x_1^p+c\cdot x^p_2\stackrel{\cR^{(0)}}{\sim}x_1^p+x^p_2$ over $\bk,$ but over $\k$ the two germs are not $\cK$-equivalent.

\item\ [The case of $j\ge1$ and  positive characteristic] Let $f(x)=x^p\sin \k[\![x]\!]$, where $char(\k)=p$.
Then $x^p\stackrel{ \cR^{(j)}_\bk}{\sim}x^p+c\cdot x^{(j+1)p}$  for any $c\sin \k$. But
 suppose the Frobenius morphism of $\k$ is non-surjective,  i.e. the field $\k$
 is non-perfect, \cite[\S26]{Matsumura}. Then $x^p\stackrel{\cR}{\not\sim}x^p+c\cdot x^{(j+1)p}$ for any $c\sin \k\smin \k^p$.

\item\ [The case of $j\ge1$ and $char=0,$ but the ring extension $\k\sset \K$ is not faithfully-flat] Let $\k$ be a ring and suppose an element $c\sin \k$ is neither invertible  nor a zero divisor.
 (E.g. $\k=\Z$ and $c=2.$) Take the extension $\k\sset \K:=\k[\frac{1}{c}] $
 and accordingly $\k[\![x]\!]\sset \K[\![x]\!].$
 Let $f(x)=c\cdot x$ and $\tf(x)=c\cdot x+x^{j+1}.$ Then   $f\stackrel{\cR^{(j)}_\K}{\sim}\tf$ for
  the filtration $(x)^\bullet\sset \K[\![x]\!],$ for $j\ge1.$ But $f\stackrel{\cR }{\not\sim}\tf.$
   Here the extension $\k\sset \K$ is not faithful, e.g. $\K\otimes\quots{\k}{(c)}=0.$
\eee

\subsection{The criterion for  equivalence with unipotent linear part}
\subsubsection{} Recall that a linear operator on a vector space, $\Phi\circlearrowright V,$ is called unipotent if the operator $\Phi-Id\circlearrowright V$ is nilpotent.
 The {\em unipotent filtration of $\Phi$} is defined as
 \beq
 V_0:=V\supset V_1:=Im(\Phi-Id)\supset V_2:=Im(\Phi-Id)^2\supset \cdots
 \eeq
Thus $\Phi$ preserves all the subspaces $V_\bullet,$ and acts as $Id$ on all the quotients $\quots{V_\bullet}{V_{\bullet+1}}.$

\subsubsection{}

Let $\k$ be an arbitrary field, and take the germs of spaces, $X=Spec(R_X),$ $Y=Spec(R_Y).$ (We assume $J_X\sseteq (x)^2$ and $J_Y\sseteq(y)^2$.)
 Take their ``first jets",
\beq
jet_1(X):=Spec(\quots{R_X}{(x)^2})\cong \k^n,\hspace{1cm}  jet_1(Y):=Spec(\quots{R_Y}{(y)^2})\cong\k^m.
\eeq
Accordingly we have the projection $\Maps\to {\rm Maps}(jet_1(X),jet_1(Y))\cong Hom_\k(\k^n,\k^m).$
 For each transformation $g\in \cG$ we take its linear part, $jet_1(g)\circlearrowright {\rm Maps}(jet_1(X),jet_1(Y)).$
  This gives the homomorphism of groups:
 \beq
 \cR\to jet_1(\cR)\le GL_\k(n),\quad
 \cL \to jet_1(\cL) \le GL_\k(m),\quad
 \cC\to jet_1(\cC)\le GL_\k(m).
 \eeq
 Similarly: $\cL\cR \to jet_1(\cL\cR) \le GL_\k(n)\times GL_\k(m)$ and
 $ \cK\to  jet_1(\cK)\le GL_\k(n)\times GL_\k(m).$
\bed
\bee
\item An element $g\in \cC$ is called unipotent if the linear operator $jet_1(g)\circlearrowright Hom_\k(\k^n,\k^m)$ is unipotent.
\hfill Take the subset of all the unipotent elements    $\cG^{unip}\sset \cG$.
\item For any field extension, $\k\sset \K,$ denote by $(\cG^{unip})_\K\sset \cG_\K$ the subset of all the unipotent elements,
 whose unipotent filtrations (on $Hom_\K(\K^n,\K^m)$) are induced from $\k,$ i.e. are of the form $\K\otimes_\k V_\bullet.$
\eee
\eed
Observe that subsets $\cG^{unip}\sset \cG$ and  $(\cG^{unip})_\K\sset \cG_\K$  are not necessarily subgroups.
 In addition  $(\cG^{unip})_\K\sseteq (\cG_\K)^{unip},$ and the inclusion is often proper.

\subsubsection{}
Let $\k\sset \K$ be an arbitrary extension of fields of char=0, e.g. $\Q\sset \C$ or $\R\sset \C.$
Let $\cG$ be one of the groups $\cR,\cL,\cL\cR,\cC,\cK.$
 In the case [$\cG=\cL,\cL\cR,$ and $f$ is analytic, and not finitely-$\cG$-determined] we assume: $f$ is of finite singularity type.
\bcor\label{Thm.Equivalence.change.base.field}
If $\tf\stackrel{(\cG^{unip})_\K}{\sim}f$ then  $\tf\stackrel{ \cG^{unip}  }{\sim}f.$
\ecor
\bpr
Suppose $\tf=g_\K\cdot f,$ for a unipotent element $g_\K\in (\cG^{unip})_\K.$
 Thus the operator $jet_1(g_\K)\circlearrowright Hom_\K(\K^n,\K^m)=:\!V$ is unipotent \wrt the filtration $\K\otimes V_\bullet,$ for some  filtration
  $(\quots{R_X}{(x)})^{\oplus m}=:\!V_0\supset V_1\supseteq\cdots\supseteq V_{m\cdot n}=0,$ by $\k$-vector spaces.
 Using the projection $R_X\stackrel{\pi}{\to}\quots{R_X}{(x)}$ we  lift  this filtration to $\RmX:$
\beq
\tM^0:=\RmX\supseteq \tM^1:=\pi^{-1}V_1\supseteq\cdots\supseteq \tM^{mn-1}:=\pi^{-1}V_{mn-1}\supseteq
\eeq
\[
\supset\tM^{mn}:=(x)\cdot \tM^0\supseteq
\cdots\supseteq \tM^{2mn-1}:=(x)\cdot \tM^{mn-1}\supset \tM^{2mn}:=(x)^2\cdot \tM^0\supseteq\cdots
\]
Thus $(x)\cdot \RmX\sseteq \tM^j\sseteq \RmX$ for $0\le j\le mn.$
For each $j\in \N$ the quotient $\quots{\tM^j}{\tM^{j+1}}$ is a vector space over $\quots{R_X}{(x)},$ of $dim=1.$
 And this filtration is equivalent to $(x)^\bullet\cdot \RmX.$

 Finally we claim: $g_\K\in \cG^{(1)}_\K$ for the filtration $\K\otimes\tM^\bullet.$ Indeed, present $g_\K$ by a power series, as in \eqref{Eq.Expanding.Autos}. The linear part $jet_1(g_\K)$ is $V_\bullet$-unipotent by the assumption. And the terms of order$\ge2$ are unipotent as well.

Finally we invoke   part 1 of Theorem \ref{Thm.Equivalence.Change.base.ring}, and the statement follows.
\epr

\vspace{-0.2cm}
\subsection{Examples and Corollaries}\label{Sec.K.vs.k.Examples.Corollaries}
\vspace{-0.4cm}

\subsubsection{Base change of families}\label{Sec.K.vs.k.Change.of.Families}

Theorem \ref{Thm.Equivalence.Change.base.ring} and Corollary \ref{Thm.Equivalence.change.base.field} address pairs of maps, $f,\tf.$  For families of maps we have a stronger result.
Let $\k=\quots{\k_o[\![t]\!]}{\ca},$ resp. $\k=\quots{\k_o\bl t\br}{\ca},$ $\k=\quots{\k_o\{ t\}}{\ca},$ where $\k_o$ is a field. Thus   elements $f\sin \Maps$ are  families  of maps over $\Spec(\k),$ denote these $f_t.$ Assume the germs $X,Y$ ``do not vary with $t$", i.e. the ideals $J_X,J_Y$ are induced from the corresponding ideals over $\k_o.$ Take any extension of fields, $\k_o\sset \K_o,$ accordingly we get:
 $\K=\quots{\K_o[\![t]\!]}{\ca},$ resp. $\K=\quots{\K_o\bl t\br}{\ca},$ $\K=\quots{\K_o\{ t\}}{\ca}.$

 In the case [$\cG=\cL,\cL\cR$ and $f$ is an analytic map] we assume: $f$ is of finite singularity type.

\bcor\label{Thm.Equivalence.for.families}
($char(\k)=0$)  If a family $f_t$ is $\cG_\K$-trivial, then it is also $\cG$-trivial.
\ecor
\bpr
Suppose the family is $\cG_\K$-trivial, i.e. $g^\K_t f_t=f_o.$ Then $g^{\K_o}_o f_o=f_o,$ where $g^{\K_o}_o\sin \cG_{\K_o}.$ Note that $\cG_{\K_o}<\cG_\K.$

Thus we replace $g^\K_t$ by $(g^{\K_o}_o)^{-1}\cdot g^\K_t.$
 Therefore we can assume $g^\K_o=Id.$

Take the filtration $(t)^\bullet\sset R_X$ and accordingly the filtration of the space of maps, induced from $\tM^\bullet:=t^\bullet\cdot \RmX.$  Then $g^\K_t\sin \cG^{(1)}_\K.$ Note also that the extension $\k \sset \K $ is faithfully-flat, see example \ref{Ex.Faithfully.Flat}. Apply part 1 of Theorem \ref{Thm.Equivalence.Change.base.ring}.
  \epr
\bex\label{Ex.Stability.base.change}
A map $X\stackrel{f}{\to}Y$ is called stable if all its local unfoldings are $\cL\cR$-trivial.
 See \cite[Chapter]{Mond-Nuno} for  $\MapC$ and \cite[\S3]{Kerner.Unfoldings} for $\MapX.$  We get (in $char=0$): $f$ is stable as a $\K$-map \iff it is stable as a $\k$-map.
\eex
\beR
This Corollary does not hold in $char=p.$ Suppose $c\in \k_o\smin \k^p_o,$ thus $\k_o$ is not a perfect field. The family $x^p+c\cdot t^p\cdot x^{dp}$ is $\cR_{\bk}$-trivial, but not $\cR_\k$-trivial.
\eeR

\subsubsection{The question `` $\k$ vs $\K$"  is finitely determined} Suppose the maps $f,\tf\sin \Maps$ are $\cG_\K$-equivalent. We prove: their $\cG$-equivalence can be verified at the level of finite jets.
  In the case [$\cG\=\cL,\cL\cR$ and $f$ is formal or analytic, and not finitely-$\cG$-determined]  we take the assumptions \eqref{Eq.SAP.Assumptions.for.L.LR.group}.

\bcor\label{Thm.K.vs.k.finitely.determined} Suppose $char(\k)=0$ and the ring extension $\k\sset \K$ is faithfully flat.
\bee
\item  Every map $f\sin \Maps$ satisfies: $\cG_\K f\!\cap \!(\{f\}+\cm^N\cdot \RmX)\sset \cG f$ for some $1\!\ll\! N\sin \N.$

\item Moreover, this number $N$ can be chosen uniformly for the whole $\cG_\K$-orbit of $f$.
\eee
\ecor
Explicitly, part 1 reads: if $f\stackrel{\cG_\K}{\sim}f+h$ for some $h\sin \cm^N\cdot \RmX$ then $f\stackrel{\cG}{\sim}f+h.$ Here $N=N_f.$

 Part 2 reads: if $f\stackrel{\cG_\K}{\sim}\tf$ for $f,\tf\sin \Maps,$ then $N_f=N_\tf.$
\bpr  Take the filtration $\cm^\bullet\cdot \RmX$ and the corresponding subgroups $\cG^{\bullet}<\cG.$
\bee
\item
Fix $1\ll N_f\sin \N$ to invoke the $\cG_\K f$-Artin-Rees lemma, \S\ref{Sec.Artin-Rees.G}. We get:
\beq
\cG_\K f\cap (\{f\}+\cm^{N_f}\cdot \RmX)\sseteq \cG^{(1)}_\K f\cap \Maps .
\eeq
 And by Theorem \ref{Thm.Equivalence.Change.base.ring} we have: $\cG^{(1)}_\K f\cap \Maps \sseteq  \cG f.$
\item
 This follows by Remark \ref{Rem.Artin.function.invariant.on.G.orbits}.
It is enough to observe (for the filtration $\cm^\bullet\cdot \RmX$): $\cG^{(o)}_\K=\cG_\K.$
\epr
\eee

\subsubsection{Finite splitting of orbits for finite field extensions $\k\sset \K$}\label{Sec.K.vs.K.Finite.Splitting.of.Orbits}
 Take a (formal, resp. analytic, Nash) map-germ $X\stackrel{f}{\to}Y,$ and a group
   $\cG=\cR,\cL,\cL\cR,\cC,\cK.$
       We get the orbit $\cG f\sset \Maps.$ Extend the field, $\k\sset \K,$ then we get the orbit $\cG_\K f.$
  Take its $\k$-part, $\cG_\K f\cap\Maps. $  It splits into $\cG$-orbits,  $\cG_\K f\cap\Maps=\amalg\cG f_\al.$
  In general the number of such orbits can  be infinite.
 Below we assume:
\bei
\item either the map $f$ is finitely-$\cG$-determined;
\item or $char=0,$ and moreover, for the case [$\cG=\cL,\cL\cR$ and the map is formal or analytic]  we take the assumptions \eqref{Eq.SAP.Assumptions.for.L.LR.group}.
\eei
\bcor\label{Thm.Finite.Splitting.orbits}
 Let the field $\k$ be  ``of type F".
Take a finite  field-extension   $\k\sset \K.$
 Then the spltting $\cG_\K f\cap \Maps=\amalg \cG f_\al$ is finite.
\ecor
A field $\k$ is called ``of type F" if $\k$ is perfect, and the Galois group $Gal(\quots{\bk}{\k})$ is ``of type F", see \cite[\S III.4.1.,pg.143]{Serre.Galois}. The simplest examples   are finite fields,  the real numbers  $ \R,$
   the field of formal power series $\k_o(\!(t)\!)$ where $\k_o=\bk_o$ is of char=0,  and p-adic fields (i.e. finite extensions of $\Q_p$).
    The field $\Q$ is not of type F.
\bpr
\bee[\bf Step 1.]
\item
Consider the $N$-jets of the spaces,
 \beq
 jet_N X:=Spec(\quots{\k[\![x]\!]}{J_X+\cm_X^{N+1}}), \hspace{1cm}
   jet_N Y:=Spec(\quots{\k[\![y]\!]}{J_Y+\cm_Y^{N+1}}),
   \eeq
   and their maps, ${\rm Maps}(jet_N X,jet_N Y).$ The latter space is a (finite-dimensional) affine variety   over $\k.$ The group-actions $\cG\circlearrowright {\rm Maps}(jet_N X,jet_N Y)$ are algebraic.
\\
We claim: it is enough to verify (for $N\gg1$) the finiteness of the decomposition
$\cG_\K jet_N f\cap {\rm Maps}(jet_N X,jet_N Y)\!= \amalg_\al \cG jet_N f_\al.$

\bei
\item This is immediate when $f$ is  finitely-$\cG$-determined.
 \item Otherwise, for $char=0,$ we invoke Corollary \ref{Thm.K.vs.k.finitely.determined}.
 \eei

\item  Thus we have the algebraic group action on a (finite-dimensional) affine variety, $\cG_\cK\circlearrowright {\rm Maps}(jet_N X_\K,jet_N Y_\K).$    It remains to recall the standard fact: (for $\k$ of $F$-type) the $\k$-part of an orbit of an algebraic group $G_\K$ (defined over $\k$) splits
   into a finite number of $\k$-orbits. See  Section III.4.4, Theorem 5 of \cite{Serre.Galois}.
\epr
\eee
\bex
Take a real map-germ $X\stackrel{f}{\to}Y.$ For [$\cG=\cL,\cL\cR$ and $f$-analytic] we assume $f$ is of finite singularity type.
  Complexifying all the objects we get  the complex orbit, $\cG_\C f\sset {\rm Maps}(X_\C,Y_\C).$ Take its real part,
   $\cG_\C f\cap \Maps.$ It splits into a finite disjont union of the real orbits.
\eex

\beR\label{Rem.Finite.Splitting.of.Orbits}
One would like to have the finite splitting  $\cG_\K f\cap \Maps=\amalg_{finite}\cG f_\al$ for more general field extensions, removing the assumptions on $\k\sset \K$ in Corollary \ref{Thm.Finite.Splitting.orbits}.
\bee[\bf i.]
\item
 For non-finite field extensions this finiteness fails trivially. Let $\k\sset \K\sseteq\bk,$ where the extension $\k\sset \K$ is not finite.
 (E.g. $\Q\sset \R$.) Fix a sequence of elements $c_\bullet\sin \K$ satisfying: $\frac{c^2_i}{c^2_j}\neq\k$ for all $i\neq j.$
   Take $f(x)=x^2.$ All the maps $c^2_j\cdot x^2$   are $\cR_\K$ equivalent to $x^2,$ but pairwise not $\cR_\k$-equivalent.
\item
When the extension $\k\sset \K$ is finite, but $\k$ is not ``of type $F$", the statement still fails, i.e.
 the $\k$-part of a $G_\K$-orbit can split into infinite number of $G_\k$-orbits.

The following example was given by Mikhail Borovoi.
 Let $\k$ be a field of $char(\k)=0,$ and assume $\sqrt{-1}\not\in \k.$ Define the action $SL(2,\k)\circlearrowright Mat^{sym}_{2\times 2}(\k)\cap SL(2,\k)$
  by  $(u,s)\to usu^T.$ Take the field extension $\K=\k(\sqrt{-1}).$
 \bei
 \item We claim:  $SL(2,\K)$ has only one orbit in $Mat^{sym}_{2\times 2}(\K)\cap SL(2,\K).$

 Indeed, the stabilizer of $\one$ is the subgroup $ SO(2,\k)<SL(2,\k) .$
This is a 1-dimensional $\k$-torus splitting over   $\K.$
In other words $SO(2,\k)\simeq \K^*,$   the multiplicative group of $\K.$

For any field extension $\k\sset L$
the set of orbits of $SL(2,L)$ in  $Mat^{sym}_{2\times 2}(L)\cap SL(2,L)$ is in a canonical bijection with the kernel
 $ker [ H^1(L,SO(2 ))\to H^1(L,SL(2)) ]$.
This  fundamental fact  is \cite{Serre.Galois}, Section I.5.4, page 50, Corollary 1 of Proposition 36.

For $L=\K$, we have $H^1(\K,SO(2 ))=H^1(\K,\K^*) = {1}$.
This is Hilbert's Theorem 90, e.g.  \cite{Serre.Fields}
  Section X.1, page 150, Proposition 2.
Thus over $\K$ we have only  one orbit, i.e.  $SL(2,\K)$ acts on $Mat^{sym}_{2\times 2}(\K)\cap SL(2,\K)$ transitively.

\item
The set of orbits of $SL(2,\k)\circlearrowright Mat^{sym}_{2\times 2}(\k)\cap SL(2,\k)$ is in a canonical bijection with $H^1(\K/\k,SO(2)) := H^1(Gal(\K/\k), SO(2,\K)),$ see \cite{Serre.Galois}, Section I.5.4, page 50, Corollary 1 of Proposition 36.
 Here  $Gal(\K/\k)$ is the group of order 2, whose non-trivial element $\tau$ acts on $SO(2,\K)\simeq\K^\times$ by sending $x\in\K^\times$ to $^\tau\!x^{-1}$, where $x\mapsto{}^\tau\! x$ is the standard action of $\tau$ on $\K$.

 As $SO(2)$  is an  abelian  group,   we can use the theory of  abelian Galois cohomology.
By definition, $H^1(Gal(\K/\k), SO(2,\k) ) = \k^\times / Nm(\K^\times),$ where the norm map acts by
 $Nm(a+b\sqrt{-1})=a^2+b^2\sin \k.$
 See \cite{Serre.Galois}, Section I.2.2, page 10.

Finally  (for $\k=\R$) the group $\R^\times/ Nm (\C^\times)$ is  of order 2.
 But (for $\k=\Q$) the group $\Q^\times/ Nm (\Q(\sqrt{-1})^\times)$ is   infinite. E.g.  let $p$ be a prime of the form $4n+3.$
  Then $p$ is not a sum of two squares in $\Z$, because   $3\ mod\ 4$ is not a sum of two squares in $\Z/4\Z.$
 Therefore $p$ is not a sum of two squares in $\Q.$
This gives infinitely many non-equivalent symmetric matrices in $Mat^{sym}_{2\times 2}(\Q).$
 \eei

\eee
\eeR

\subsubsection{A global implication}\label{Sec.K.vs.k.Global}
Below we work with graded rings, $R_X=\quots{\k[x]}{J_X},$ $R_Y=\quots{\k[y]}{J_Y},$ where $J_X,J_Y$ are generated by homogeneous polynomials.
 Take a non-zero graded map $X\stackrel{f}{\to}Y.$ Projectivizing all the objects we get the (rational) map of projective schemes $\P X\stackrel{\P f}{\to }\P Y,$ see e.g. pg.160-169 of \cite{Hartshorne}.

  Let $\cG$ be one of the groups $\cR,\cL,\cL\cR.$ Recall the standard fact: if two graded maps  $X\stackrel{f,\tf}{\to}Y$ are formally-$\cG$-equivalent, then they are graded-$\cG$-equivalent.   Therefore the maps $\P f,\P\tf$ are projectively equivalent.

Part 1 of Theorem \ref{Thm.Equivalence.Change.base.ring} is trivial  in this case. Indeed, a graded automorphism that is unipotent (for the graded filtration) is necessarily identity. Part 2   is well-known.

But Corollary \ref{Thm.Equivalence.for.families} gives the (well-known?) statement:
\bcor
(char=0) Take a family of projective morphisms $\P X\stackrel{\P f_t}{\to }\P Y.$ If this family is globally-$\K$-trivial, then it is globally-$\k$-trivial.
\ecor

\end{document}